\documentclass[11pt,reqno]{amsart}
\usepackage {amssymb}
\usepackage {amsmath}
\usepackage {bbm}
\usepackage{amsthm}
\usepackage{mathtools}
\usepackage{graphicx}
\usepackage {amscd}
\usepackage{makecell}
\usepackage {epic}
\usepackage{array}
\usepackage{enumerate}
\usepackage{mathrsfs}
\usepackage{subcaption} 
\usepackage{etaremune}
\usepackage{stmaryrd}

\DeclareFontFamily{U}{dutchcal}{\hyphenchar\font=-1}
\DeclareFontShape{U}{dutchcal}{m}{n}{ <-> dutchcal-r}{}
\DeclareSymbolFont{dutchletters}{U}{dutchcal}{m}{n}

\DeclareMathSymbol{\mathdutchcal}{0}{dutchletters}{"41} 
\newcommand{\dcal}[1]{\text{\usefont{U}{dutchcal}{m}{n}#1}}

\usepackage[dvipsnames]{xcolor}
\usepackage[colorlinks,citecolor=CornflowerBlue,linkcolor=RoyalPurple,urlcolor=Plum,pagebackref]{hyperref}
\usepackage{cleveref}

\usepackage{enumerate}
\usepackage{tikz}
\usepackage{tikz-cd}
\usetikzlibrary{arrows.meta}
\usepackage{verbatim,color,geometry}
\usepackage[all]{xy}
\usepackage{bm}
\usepackage{longtable}
\usepackage{aliascnt}
\usetikzlibrary{calc}
\usepackage{textcomp}
\usepackage{tablefootnote}
\usepackage{float}

\newcommand{\mtf}[1]{\mathfrak{#1}}

\newcommand{\wh}[1]{\widehat{#1}}
\newcommand{\ove}[1]{\overline{#1}}

\newcommand{\cC}{\dcal{C}}
\newcommand{\cD}{\dcal{D}}

\newcommand{\cL}{\dcal{L}}

\newcommand{\cO}{\dcal{O}}
\newcommand{\cP}{\dcal{P}}

\newcommand{\cR}{\dcal{R}}

\newcommand{\cU}{\dcal{U}}

\newcommand{\cZ}{\dcal{Z}}
\newcommand{\cX}{\dcal{X}}
\newcommand{\cY}{\dcal{Y}}

\newcommand\MM{{\mathcal{M}}}

\newcommand{\bP}{\mathbb{P}}
\newcommand{\bC}{\mathbb{C}}

\newcommand{\bA}{\mathbb{A}}
\newcommand{\bQ}{\mathbb{Q}}
\newcommand{\bZ}{\mathbb{Z}}

\newcommand{\bG}{\mathbb{G}}
\newcommand{\bF}{\mathbb{F}}

\newcommand{\bN}{\mathbb{N}}

\newcommand{\klt}{\mathrm{klt}}
\newcommand{\Bl}{\mathrm{Bl}}

\newcommand{\KSBA}{\mathrm{KSBA}}

\newcommand{\lc}{\mathrm{lc}}

\newcommand{\Sym}{\mathrm{Sym}}

\newcommand{\sX}{\mathscr{X}}
\newcommand{\sY}{\mathscr{Y}}
\newcommand{\sZ}{\mathscr{Z}}

\DeclareMathOperator{\Aut}{Aut}

\DeclareMathOperator{\vol}{vol}

\DeclareMathOperator{\Fut}{Fut}

\DeclareMathOperator{\lct}{lct}

\DeclareMathOperator{\mult}{mult}

\DeclareMathOperator{\GL}{GL}

\DeclareMathOperator{\bb}{BB}

\DeclareMathOperator{\coeff}{coeff}

\DeclareMathOperator{\ord}{ord}

\DeclareMathOperator{\Proj}{Proj}

\DeclareMathAlphabet{\mathbbb}{U}{bbold}{m}{n}

\newcommand{\mb}[1]{\mathbb{#1}}

\newcommand{\SL}{\mathrm{SL}}

\newcommand{\Chow}{\mathrm{Chow}}

\newcommand{\bT}{\mathbb{T}}

\newcommand{\tQ}{\widetilde{Q}}

\newcommand{\ofM}{\overline{\mathfrak{M}}}

\newcommand{\GIT}{\mathrm{GIT}}

\newcommand{\fM}{\mathfrak{M}}
\newcommand{\fF}{\mathfrak{F}}

\newcommand{\K}{\mathrm{K}}

\newcommand{\Def}{\mathrm{Def}}
\newcommand{\sst}{\mathrm{ss}}

\numberwithin{equation}{section}

\newtheorem{prop}{Proposition}[section]

\newcommand{\newaliastheorem}[3]{%
  \newaliascnt{#1}{prop}%
  \newtheorem{#1}[#1]{#2}%
  \aliascntresetthe{#1}%
  \crefname{#1}{#2}{#3}%
  \Crefname{#1}{#2}{#3}%
}

\newaliastheorem{thm}{Theorem}{Theorems}
\newaliastheorem{theorem}{Theorem}{Theorems}
\newaliastheorem{lem}{Lemma}{Lemmas}
\newaliastheorem{lemma}{Lemma}{Lemmas}
\newaliastheorem{cor}{Corollary}{Corollaries}
\newaliastheorem{corollary}{Corollary}{Corollaries}
\newaliastheorem{theodef}{Theorem-Definition}{Theorem-Definitions}
\newaliastheorem{prop-def}{Proposition-Definition}{Proposition-Definitions}
\newaliastheorem{convention}{Convention}{Conventions}
\newaliastheorem{conj}{Conjecture}{Conjectures}
\newaliastheorem{conjecture}{Conjecture}{Conjectures}

\theoremstyle{definition}
\newaliastheorem{defi}{Definition}{Definitions}
\newaliastheorem{defn}{Definition}{Definitions}
\newaliastheorem{exa}{Example}{Examples}
\newaliastheorem{exam}{Example}{Examples}
\newaliastheorem{expl}{Example}{Examples}
\newaliastheorem{rem}{Remark}{Remarks}
\newaliastheorem{rmk}{Remark}{Remarks}
\newaliastheorem{remark}{Remark}{Remarks}
\newaliastheorem{que}{Question}{Questions}

\crefname{prop}{Proposition}{Propositions}
\Crefname{prop}{Proposition}{Propositions}

\makeatletter

\newcommand*{\@rowstyle}{}

\newcommand*{\rowstyle}[1]{%
  \gdef\@rowstyle{#1}%
  \@rowstyle\ignorespaces%
}

\newcolumntype{=}{%
  >{\gdef\@rowstyle{}}%
}

\newcolumntype{+}{%
  >{\@rowstyle}%
}

\makeatother

\DeclareMathOperator{\CY}{CY}

\DeclareMathOperator{\Hodge}{Hodge}

\newcommand{\sheafHom}{\mathscr{H}\text{\kern -3pt {\calligra\large om}}\,}

\DeclareMathOperator{\divisor}{\mathrm{div}}

\newcommand{\fN}{{\mathfrak{N}}}

\DeclareFontFamily{U}{dutchcal}{\hyphenchar\font=-1}
\DeclareFontShape{U}{dutchcal}{m}{n}{ <-> dutchcal-r}{}
\DeclareSymbolFont{dutchletters}{U}{dutchcal}{m}{n}

\DeclareMathSymbol{\mathdutchcal}{0}{dutchletters}{"41}

\theoremstyle{definition}

\theoremstyle{remark}

\title{Birational Geometry of K-Moduli Spaces of Genus Four Curves}
\date{\today}

\author{Kenneth Ascher}
\address{340 Rowland Hall, Irvine, CA 92697, USA}
\email{kascher@uci.edu}

\author{Kristin DeVleming}
\address{2985 Muir Ln, La Jolla, CA 92093, USA}
\email{kedevleming@ucsd.edu}

\author{Hongyu Luo}
\address{Gordon Palmer Hall, 505 Hackberry Lane, Tuscaloosa, AL 35401, USA}
\email{hluo13@crimson.ua.edu}

\author{Yuchen Liu}
\address{Lunt Hall, 2033 Sheridan Road,
Evanston, IL 60208, USA}
\email{yuchenl@northwestern.edu}

\author{Junyan Zhao}
\address{William E. Kirwan Hall, 4176 Campus Dr, College Park, MD 20742, USA}
\email{jzhao81@umd.edu}

\begin{document}

\begin{abstract}
We study the wall-crossings for K-moduli spaces of $(3,3)$-curves on $\bP^1\times\bP^1$. Together with the boundary polarized Calabi--Yau moduli space, these spaces interpolate between the GIT moduli space of $(3,3)$-curves and the Baily--Borel compactification of the moduli space of K3 surfaces arising as cyclic triple covers of $\bP^1\times\bP^1$ branched along such curves. As applications, we describe the corresponding birational models in the Hassett--Keel--Looijenga program and realize the K-moduli space of degree one del Pezzo surfaces through explicit wall-crossings.
\end{abstract}

\maketitle

\tableofcontents

\section{Introduction}

Constructing compactifications of moduli spaces is a fundamental problem in algebraic geometry. Different compactifications of a given moduli space of varieties parametrize different singular objects on the boundary, often reflecting distinct geometric features of the varieties in question. For smooth curves of genus $g\geq2$, the classical Deligne--Mumford compactification $\ove{M}_g$ adds stable nodal curves to yield a proper moduli space. In low genus, concrete projective models of curves provide alternative compactifications through geometric invariant theory and moduli of pairs.  Given two different compactifications arising from different stability conditions, it is then natural to ask for a geometric description of the birational map between them. 

In this paper, we apply these ideas to compactifications of the moduli space of smooth $(3,3)$ curves on the quadric surface $\bP^1 \times \bP^1$.  Observe that a general canonical genus four curve is a $(2,3)$-complete intersection in $\bP^3$ lying on a smooth quadric, and hence can be viewed as a $(3,3)$-curve on $\bP^1\times\bP^1$.  This observation leads us to consider various compactifications of geometric significance:
\begin{itemize}
\item the GIT and VGIT quotients of $(3,3)$-curves and $(2,3)$-    complete intersections;
    \item the K-moduli, KSBA, and boundary polarized Calabi Yau compactifications of pairs of the form $(\bP^1\times\bP^1,cC)$ for $0<c\leq 1$;
    \item the Baily--Borel and Looijenga compactifications of Kond\={o}'s ball quotient model of K3 surfaces arising as cyclic covers of $\bP^1\times\bP^1$ branched along $(3,3)$ curves; and
    \item some Hassett--Keel birational models of $\overline{M}_4$.
\end{itemize}

Several of these constructions and their relationships have been studied previously; see \Cref{sec:Prior-work}. Our main approach in this paper is to consider the K-moduli compactifications and their wall-crossings as $c\in(0,\frac{2}{3})$ varies.  In particular, we will prove that the K-moduli spaces interpolate between the GIT compactification and the Baily--Borel compactification.

Let $\ove{\fM}^{\GIT}\coloneqq
    |\cO_{\bP^1\times\bP^1}(3,3)|
    \sslash \Aut(\bP^1\times\bP^1)$ be the GIT moduli space of $(3,3)$-curves. For $c\in(0,\tfrac23)\cap\bQ$, let $\MM_4^\K(c)$ be the moduli stack of K-semistable log Fano pairs admitting $\bQ$-Gorenstein smoothings to $(\bP^1\times\bP^1,cC)$, where $C$ is a $(3,3)$-curve, and let $\ove{\fM}_4^\K(c)$ be its good moduli space; see \Cref{def: k moduli spaces of bidegree 3}. For each K-moduli wall $c_i$, let $E_i^{\pm}\subseteq \ove{\fM}_4^\K(c_i\pm\epsilon\big)$ denote the local exceptional loci of the wall-crossing morphisms at $c_i$. See \Cref{thm:wall-crossing-structure} and the subsequent discussion for the relevant terminology. Roughly, the reader should think of $E_i^-$ as `the locus of curves that get replaced' when crossing the wall, and $E_i^+$ as `the locus of replacements.' Our first theorem describes all of the K-moduli spaces, their wall-crossings as $c$ varies, and the relationship with the aforementioned GIT quotient.

\begin{theorem}\label{Main Theorem 1.1} For any $c\in (0,\tfrac{2}{3})$, the K-moduli stack $\MM_4^\K(c)$ is smooth, and the K-moduli space $\ofM_4^{\K}(c)$ is normal. Moreover, 
\begin{enumerate}
   \item[\textup{(1)}] for any $c\in (0,\frac{1}{6})\cap \bQ$, there is an isomorphism $\ove{\fM}_4^\K(c)\cong \ove{\fM}^\GIT$;
   \item[\textup{(2)}] there are exactly six K-moduli walls \[
\big(c_1,c_2,\cdots, c_6\big)\ =\  \left( \frac{1}{6}, \ \frac{1}{3},\  \frac{16}{33}, \ \frac{1}{2}, \ \frac{11}{21}, \ \frac{17}{30} \right)\] in $(0,\tfrac{2}{3})$ with exactly seven wall-crossing centers, among which $c_1,c_3$ are divisorial contractions and the others are flips;
   \item[\textup{(3)}] for any $c \in (0,\tfrac{2}{3})$ and any $[(X,cD)]\in \MM_4^\K(c)$, the surface $X$ is isomorphic to one of the following: $\bP^{1} \times \bP^{1}$, $\bP(1,1,2)$, $\bP(1,2,9)$ or its partial smoothing $X_{10}$.
\end{enumerate}
The geometry of each K-moduli wall-crossing is summarized in \Cref{table:K-moduli walls}.
\end{theorem}

\begin{table}[htbp]
\centering
\renewcommand{\arraystretch}{1.2}
\setlength{\tabcolsep}{5pt}

\begin{tabular}{
|c|c|c|c|c|c|
}
\hline
$i$ & $c_i$
& $\dim E_i^-$ 
& $E_i^-$ 
& $E_i^+$ & $\dim E_i^+$ \\
\hline

1 & $\frac{1}{6}$
  & 0
  & triple conic 
  & $(\bP(1,1,2), V(s^3 + f_4(u,v)s + f_6(u,v) ))$ & 8\\\hline

2 & $\frac{1}{3}$
  & 1
  & double conics  
  & $(\bP(1,1,2), V(uvs^2 + (au^4+bv^4)s + f_6(u,v)))$ & 7\\\hline

3 & $\frac{16}{33}$
  & 0
  & $A_9^{\rm sep}$-curve  
  & $(X_{10}, V(w - g(x,y)))$ & 8 \\\hline

4 & $\frac{1}{2}$
  & 1
  & $A_8$-curves  
  & $(\bP(1,2,9), V(s^2 + v^9 + u^4g(u,v) ))$& 7 \\\hline

4' & $\frac{1}{2}$
  & 2
  & $A_7^{\rm sep}$-curves  
  & $(\bP(1,1,2), V(u^2s^2 + v^4s + g(u,v) ))$& 6 \\\hline

5 & $\frac{11}{21}$
  & 2
  &  $A_7^{\rm nsep}$-curves  
  & $(\bP(1,2,9), V(s^2 + u^2v^8 + u^6g(u,v) ))$& 6 \\\hline

6 & $\frac{17}{30}$
  & 3
  & $A_6$-curves  
  & $(\bP(1,2,9), V(s^2 + u^4v^7 + u^8g(u,v) ))$ & 5\\

\hline
\end{tabular}

\caption{Geometry of wall-crossings for $\ove{\fM}^\K_4(c)$}
\label{table:K-moduli walls}
\end{table}

To connect these spaces to the Baily--Borel compactification $\ove{\mtf{F}}^{\bb}$, we use the moduli space $\ove{\fM}_4^{\CY}$ of boundary polarized Calabi--Yau pairs introduced in \cite{ABB+23}, which parametrizes slc degenerations of $(\bP^1\times\bP^1,\tfrac23C)$; see \Cref{def:CY4stack}. By \cite[Theorem 6.16]{ADLW25}, $\ove{\fM}_4^{\CY}\simeq\ove{\mtf{F}}^{\bb}$. This allows us to view the Baily--Borel compactification as a moduli space of pairs, and then the K-moduli wall-crossings, together with the Calabi--Yau limit at $c=\tfrac23$, connect the GIT quotient to the Baily--Borel compactification.

More generally, we can also identify the K-moduli spaces with birational models of $\ove{\mtf{F}}^{\bb}$ arising in the \emph{Hassett--Keel--Looijenga program}. Following \cite{CMJL2012}, let $H_n$, $H_h$, and $H_v$ denote the nodal, hyperelliptic, and Gieseker-Petri divisors, respectively, and let $\lambda_{\Hodge}$ be the Hodge line bundle.

\begin{theorem}\label{thm:intro-K3s}
    Let $\ove{\fM}_4^\K(c)$ denote the good moduli space parametrizing K-polystable log Fano pairs $(X, cD)$ admitting a $\mathbb{Q}$-Gorenstein smoothing to $(\mathbb{P}^{1} \times \mathbb{P}^{1},cC)$, where $C$ is a $(3,3)$-curve. 
    \begin{enumerate}
        \item[\textup{(1)}] On $\ove{\mtf{F}}^{\bb}$, one has the relation \begin{equation}
    102\lambda_{\Hodge} = 2H_n + 9H_v + 28H_h.
\end{equation}
\item[\textup{(2)}] There is an isomorphism \begin{equation*}
    \ove{\fM}_4^{\K}(c) \simeq \Proj R\big(\mtf{F},\ (1-\tfrac{3}{2}c)(4H_{v} + \tfrac{32}{3}H_{h}) +36c\lambda_{\Hodge}\big).
\end{equation*}
Moreover, the CM $\bQ$-line bundle $\Lambda_c$ on $\ove{\fM}_4^{\K}(c)$ is proportional to\\ $(p_c^{-1})_{*}\big((1-\tfrac{3}{2}c)(4H_{v} + \tfrac{32}{3}H_{h}) +36c\lambda_{\Hodge}\big)$, where $p_c$ is the birational period map $\ove{\fM}_4^{\K}(c) \dashrightarrow \mtf{F}$.
\item[\textup{(3)}] For $0 < \epsilon \ll 1$, the moduli space $\ove{\fM}_4^{\CY}\simeq\ove{\mtf{F}}^{\bb}$ is the ample model of the Hodge line bundle on $\ove{\fM}_4^{\K}(\frac{2}{3}-\epsilon)$.
    \end{enumerate}
For a detailed description of (3), see Table \ref{table:K-to-CY-wall-crossing}.
\end{theorem}

As a further application, we recover the K-moduli space of degree one del Pezzo surfaces $\fM^{\K}_{\textup{dP}^1}$ as considered by Odaka--Spotti--Sun through explicit wall-crossings. This explicitly realizes the two-step birational modification described in \cite{OSS16}: a blowup followed by a small contraction.

\begin{theorem}\label{thm:intro-dp1}
Let ${\fN^{\GIT}}$ be the GIT moduli space of degree six curves on $\bP(1,1,2)$ that do not pass through the unique singular point. The K-moduli space $\fM^{\K}_{\textup{dP}^1}$ of del Pezzo surfaces of degree one fits into the following diagram: \[\begin{tikzcd}[ampersand replacement=\&]
	\& {\wh{\fN}} \&\& \\
	{\ove{\fN}^{\GIT}} \&\& {\ove{\fN}^{\K}(\frac{1}{2})} \& {\ove{\fM}^{\K}_{\textup{dP}^1}}
	\arrow["\phi"', from=1-2, to=2-1]
	\arrow["\psi", from=1-2, to=2-3]
	\arrow["\iota", from=2-3, to=2-4]
\end{tikzcd}\]
   \begin{itemize}
       \item $\phi$ is a divisorial contraction, replacing the union of a double conic with another transverse conic by curves with a node at the vertex of $\bP(1,1,2)$;
       \item $\psi$ is a small contraction, which 
       \begin{enumerate}
           \item[\textup{(1)}] contracts the rational curve parametrizing curves with a separating $A_7$-singularity to the point representing the curve  $\big(s(su^2+v^4) = 0\big)$ on $\bP(1,1,2)_{[u:v:s]}$;
           \item[\textup{(2)}] is a local isomorphism near the point representing the curve with an $A_8$-singularity, sending it to the curve $(s^2+v^9 = 0)$ on $\bP(1,2,9)_{[u:v:s]}$.
       \end{enumerate}
       \item $\iota$ is an isomorphism induced by taking the double cover branched over the given curve.
   \end{itemize}
\end{theorem}

\subsection{Prior and related work}\label{sec:Prior-work}
The wall-crossing theory for K-moduli of log Fano pairs was developed in \cite{ADL24,Zho23,LZ26}, which also determines all K-moduli spaces of plane curves of degrees $4$, $5$, and $6$. In degrees $4$ and $6$, these spaces interpolate between GIT quotients and Baily--Borel compactifications of the associated K3 moduli spaces. Our work extends this picture to $(3,3)$-curves.

Comparisons between K-moduli and (V)GIT have been established in many settings; see, for example, \cite{MM20,OSS16,SS17,Li19,Fuj18,GMG21,ADL21,zha22,zha23b,LZ24,zhao_2023,MPZ24,CKSZ26,Zha26,PSW23,ADLW25}. Of particular relevance is \cite{ADL21}, which identifies the K-moduli spaces of $(4,4)$-curves with VGIT quotients of $(2,4)$-complete intersections and, together with work of Laza--O'Grady, relates them to the Baily--Borel compactification of the moduli space of hyperelliptic quartic K3 surfaces. However, for $(3,3)$-curves, the VGIT quotients of \cite{Fed12,CMJL14} do not always coincide with the K-moduli spaces: the wall at $c=\frac{16}{33}$ introduces curves on $X_{10}$ that are not $(2,3)$-complete intersections.

Our work is also closely related to \cite{ADLW25}, which completes the Hassett--Keel program for $\ove{M}_4$. The K-moduli spaces studied here do not, in general, coincide with the Hassett--Keel models. In particular, the wall at $c=\frac{16}{33}$ extracts a divisor birational to a $\bP^1$-fibration over the hyperelliptic locus; this divisor does not occur in the Hassett--Keel models. However, by \cite{Fed12, CMJL14, ADL21}, for $c \in (0, \frac{16}{33}) \cap \mathbb{Q}$, the K-moduli spaces coincide with the last Hassett-Keel birational models of $\ove{M}_4$.   A birational transform of the same divisor appears in the KSBA compactification studied by Deopurkar--Han \cite{DH21}. The two pictures are linked by the boundary polarized Calabi--Yau moduli space, identified with Kond\={o}'s Baily--Borel compactification \cite{Kon02}. In particular, the boundary polarized Calabi--Yau moduli space connects our final K-moduli space with the KSBA compactification of \cite{DH21}.

\subsection{Outline of the paper}
\S\ref{sec:prelim} reviews the relevant background and identifies the surfaces appearing in our K-moduli spaces. \S\ref{sec: first wall} establishes the comparison with GIT for small coefficients and describes the first wall-crossing. \S\ref{sec:walls of P(1,1,2)}, \S\ref{sec:walls of X10}, and \S\ref{sec: walls of P(1,2,9)} determine the remaining walls associated with $\bP(1,1,2)$, $X_{10}$, and $\bP(1,2,9)$, respectively, and complete the proof of \Cref{Main Theorem 1.1}. Finally, \S\ref{Calabi-Yau section} applies these results to moduli of K3 surfaces and degree one del Pezzo surfaces.

\subsection{Conventions and notations}

We adopt the following conventions throughout this paper.
\begin{itemize}
    \item We work over the field $\bC$ of complex numbers.
    \item We follow the conventions of \cite{KM98} regarding singularities of varieties and log pairs.
    \item We use $([x:y],[z:w])$ to denote the homogeneous coordinates on $\bP^1\times\bP^1$, use $[x:y:z:w]$ to denote the weighted homogeneous coordinates on $\bP(1,1,5,9)$, and use $[u:v:s]$ to denote the weighted homogeneous coordinates on $\bP(1,1,2)$ and $\bP(1,2,9)$.
\end{itemize}

For the reader’s convenience, we collect the notation for the moduli stacks and spaces most frequently used throughout the paper.

\renewcommand{\arraystretch}{1.5}
\begin{center}
\begin{longtable}{|c|l|}
    \hline \textbf{Notation} & \textbf{Definition/Description}   \\ 
    \hline $\MM^\K_{4}(c) $&   K-moduli stack of pairs deforming to $(\bP^1\times\bP^1,cC)$ for $C\in |\cO_{\bP^1\times\bP^1}(3,3)|$\\
    \hline $\ove{\fM}^\K_{4}(c) $&   Good moduli space of $\MM^\K_{4}(c) $ \\  \hline $\ove{\fM}^\K_{4} $&   K-moduli space $\ove{\fM}^\K_{4}(\tfrac{2}{3}-\epsilon)$ for $0<\epsilon\ll1$  \\ \hline
    $\ove{\fM}^{\GIT}$&   GIT quotient of $|\cO_{\bP^1\times\bP^1}(3,3)|$ by $\SL_2 \times \SL_2 \rtimes \bZ_2$-action\\\hline
    $\ove{\fM}^{\GIT}(t)$&   GIT moduli space of $(2,3)$-c.i. in $\bP^3$ with linearization of slope $t$  \\\hline
    $\ove{\fM}^{\CY}_4$&  Good moduli space of log Calabi--Yau pairs deforming to $(\bP^1\times\bP^1,\frac{2}{3}C)$ \\\hline
    $\fF$ & Moduli space of K3 surfaces arising as cyclic triple cover of $\bP^1\times\bP^1$
    \\\hline

    $\ove{\fF}^{\bb}$ &   Kond\=o's Baily--Borel  compactification of $\fF$
    \\\hline
$\phi_{c_i}^{\pm}$ &  Wall-crossing morphism $\phi_{c_i}^{\pm}:\ove{\fM}^{\K}_4(c_i\pm\epsilon)\to \ove{\fM}^{\K}_4(c_i)$
    \\\hline
    $E_i^{\pm}$ &  The $\phi_{c_i}^{\pm}$-exceptional loci $E_i^{\pm}\subseteq \ove{\fM}^{\K}_4(c_i\pm\epsilon)$
    \\\hline
    $\ove{\fM}^{\K}_{\textup{dP}^1}$ &  K-moduli space of del Pezzo surfaces of degree one.
    \\\hline
    
 \caption{Notations and descriptions of moduli stacks}
    \label{tab:notations}
\end{longtable}

\end{center}

\subsection{Acknowledgments}
We would like to thank I. Cheltsov and P. Engel for useful discussions. KA, KD, and YL began work on this project at the American Institute of Mathematics as part of the AIM SQuaREs program. KA was partially supported by NSF grants DMS-2140781 and DMS-2302550. KD was partially supported by NSF grant DMS-2550445 and NSF CAREER grant DMS-2542575. YL was partially supported by NSF CAREER Grant DMS-2237139 and an AT\&T Research Fellowship from Northwestern University.

\subsection{AI usage disclosure} 

GenAI was only used for editorial assistance in preparing and auditing this manuscript. All mathematical content was developed by the authors, who take full responsibility for the content and correctness of the paper.

\section{Preliminaries}\label{sec:prelim}

\subsection{K-stability and K-moduli theory}

We begin by recalling the basic concepts and results concerning the K-stability for log Fano pairs.

\begin{defn} 
    Let $(X,D)$ be a projective klt pair. We call $(X,D)$ a \emph{log Fano pair} if $-(K_X+D)$ is ample. A projective variety $X$ is called a \emph{$\bQ$-Fano variety} if $(X,0)$ is a log Fano pair.
\end{defn}

\begin{defn} 
Let $(X,D)$ be an $n$-dimensional log Fano pair, and $E$ be a prime divisor over $X$, i.e. a divisor on a normal projective variety $Y$ admitting a birational morphism $\pi:Y\rightarrow X$. The \emph{log discrepancy} of $(X,D)$ with respect to $E$ is $$A_{X,D}(E) \ \coloneqq \ 1+\coeff_{E}\big(K_Y+\pi_*^{-1}D-\pi^{*}(K_X+D)\big),$$ the \emph{S-invariant} of $(X,D)$ with respect to $E$ is $$S_{X,D}(E) \ \coloneqq \ \frac{1}{(-K_X-D)^n}\int_{0}^{\infty}\vol_Y\big(-\pi^{*}(K_X+D)-tE\big)\,dt,$$ and the \emph{$\beta$-invariant} of $(X,D)$ with respect to $E$ is $$\beta_{X,D}(E)\ \coloneqq \ A_{X,D}(E)-S_{X,D}(E).$$
\end{defn}

\begin{theodef} \label{valuative criterion} \textup{(cf.\ \cite{Fuj19,Li17,BX19})} A log Fano pair $(X,D)$ is 
\begin{enumerate}
    \item[\textup{(1)}] \textup{K-semistable} if and only if $\beta_{X,D}(E)\geq 0$ for every prime divisor $E$ over $X$;
    \item[\textup{(2)}] \textup{K-stable} if and only if $\beta_{X,D}(E)>0$ for every prime divisor $E$ over $X$;
    \item[\textup{(3)}] \textup{K-polystable} if and only if it is K-semistable and every $\mb{G}_m$-equivariant K-semistable degeneration of $(X,D)$ is isomorphic to $(X,D)$ itself.
    \item[\textup{(4)}] \textup{K-unstable} if it is not K-semistable.
\end{enumerate}
\end{theodef}

\begin{defn}  
Let $(X,D)$ be a log Fano pair, and let $\bT$ be a maximal torus in $\Aut(X,D)$. The \textit{complexity} of $(X,D)$ is the number $\dim X - \dim\bT$. A $\bT$-invariant prime divisor $E$ on $X$ is called \textit{vertical} if a maximal $\bT$-orbit in $E$ has the same dimension as the torus $\bT$; otherwise, $E$ is called \textit{horizontal}.
\end{defn}

When $(X,D)$ is a pair of complexity 1, the following result characterizes its K-polystability.  

\begin{theorem}[\textup{cf. \cite[Theorem 1.31]{ACC23}}]\label{thm:complexity-one} 
Let $(X,D)$ be a log Fano pair of complexity $1$. Then $(X,D)$ is K-polystable if and only if the following conditions hold:
\begin{enumerate}
    \item[\textup{(1)}] $\beta_{X,D}(F) > 0$ for every vertical $\bT$-invariant prime divisor $F$ on the variety $X$;
    \item[\textup{(2)}] $\beta_{X,D}(F) \geq 0$ for every horizontal $\bT$-invariant prime divisor $F$ on the variety $X$;
    \item[\textup{(3)}] The $\beta$-invariant vanishes for every valuation induced by a one-parameter subgroup of $\bT$.
\end{enumerate}
\end{theorem}

\begin{prop}[\textup{cf. \cite[Proposition 2.13 ]{ADL24}}]\label{Interpolation of K-stability ADL}
Let $X$ be a $\bQ$-Fano variety, and let $D$ and $\Delta$ be two effective $\bQ$-divisors on $X$ satisfying
the following properties:
\begin{itemize}
    \item Both $D$ and $\Delta$ are $\bQ$-linearly equivalent to a multiple of $-K_X$;
    \item $-K_X - D$ is ample, and $-K_X - \Delta$ is nef;
    \item $(X,D)$ and $(X,\Delta)$ are K-(poly/semi)stable and K-semistable log Fano pairs, respectively.
\end{itemize}
Then the following holds:
\begin{enumerate}
    \item[\textup{(1)}] If $D \neq 0$, then $(X,tD + (1-t)\Delta)$ is K-(poly/semi)stable for any $t \in (0,1]$;
    \item[\textup{(2)}] if $D = 0$, then $(X,(1-t)\Delta)$ is K-semistable for any $t \in (0,1]$;
    \item[\textup{(3)}] if $\Delta \sim_{\mathbb{Q}} -K_X$ and $(X,\Delta)$ is klt, then $(X,tD+ (1-t)\Delta)$ is K-stable for any $ t \in (0,1)$.
\end{enumerate}
\end{prop}

Now we introduce the K-moduli theory; we refer the reader to \cite[Chapter 7 \& 8]{Xu24} for more details.

\begin{theorem}[K-moduli Theorem]\label{thm:Kmodthm} Let $c, r \in \mathbb{Q}_{> 0}$ with $c < \min\{1, r^{-1}\}$, and $\chi$ be the Hilbert polynomial of an anti-canonically polarized smooth Fano variety. Then the moduli pseudo-functor $\MM^\K_{\chi,r}(c)$ sending a reduced scheme $S$ to 
\[
\left\{
(\mathcal{X}, D)/S \;\middle|\;
\begin{array}{l}
(\mathcal{X}, cD)/S \text{ is a } \mathbb{Q}\text{-Gorenstein smoothable  family of log Fano pairs},\\
D \sim_{S,\mathbb{Q}} -rK_{\mathcal{X}/S}, \text{ each fiber } (\mathcal{X}_s,cD_s) \text{ is K-semistable,}\\
\text{and } \chi(\mathcal{X}_s,\cO_{\mathcal{X}_s}(-mK_{\mathcal{X}_s}))=\chi(m)
\text{ for } m \text{ sufficiently divisible}.
\end{array}
\right\}.
\]
is represented by a reduced Artin stack $\mathcal{M}^\K(c)$ of finite type over $\bC$, whose points parametrize K-semistable $\bQ$-Gorenstein smoothable log Fano pairs. Moreover, $\mathcal{M}^\K(c)$ admits a good moduli space $\overline{\fM}^\K(c)$, which is a reduced projective scheme, parameterizing the K-polystable log Fano pairs. 
\end{theorem}

The definition of a family of pairs is subtle, even over a reduced base scheme; we refer the reader to \cite[Chapter 4]{Kol23}.

The following structure theorem describes how the K-moduli stacks and spaces vary with $c$.

    \begin{theorem}[\textup{\cite[Theorem 1.2]{ADL24}}]\label{thm:wall-crossing-structure}
    There exist rational numbers
\[
0 \ =\  c_0 \ <\  c_1 \ < \  c_2 \ < \ \cdots \ <\  c_n \ =\  \min \{1, r^{-1} \}
\]
such that for every \(0 \le j < n\), the K-moduli stacks
\(\mathcal{M}^\K(c)\) are independent of the choice of
\(c \in (c_j, c_{j+1})\). Moreover, for every \(0 < j < n\)
and \(0 < \epsilon \ll 1\), there are open immersions
\[
\mathcal{M}^\K(c_j - \epsilon)
\ \hookrightarrow\ 
\mathcal{M}^\K(c_j)
\ \hookleftarrow\ 
\mathcal{M}^\K(c_j + \epsilon),
\]
which descend to projective morphisms
\[
\overline{\fM}^\K(c_j - \epsilon)
\ \longrightarrow\ 
\overline{\fM}^\K(c_j)
\ \longleftarrow\ 
\overline{\fM}^\K(c_j + \epsilon).
\]
Moreover, each of these wall-crossing morphisms admits a local VGIT presentation as in \cite[(1.2)]{AFS173}.
\end{theorem}

Suppose that $(X,c_0D)$ is a $\bQ$-Gorenstein smoothable K-polystable log Fano pair with $D\sim_{\bQ}-rK_X$ and $c_0\in(0,\min\{1,r^{-1}\})\cap\bQ$. By \cite[Proposition 3.18]{ADL24}, the set
\[
\{\,c\in(0,\min\{1,r^{-1}\}) : (X,cD)\text{ is K-polystable}\,\}
\]
is either the singleton $\{c_0\}$ or an open interval containing $c_0$. The critical values $c_i$ in the preceding theorem are the endpoints of these intervals; equivalently, they are the values for which there exists a pair $(X,D)$ such that $(X,cD)$ is K-polystable precisely when $c=c_i$. We call each such $c_i$ a \textit{K-moduli wall} (or simply a \textit{wall}), and the intervals between consecutive walls \textit{K-moduli chambers} (or simply \textit{chambers}). Within each chamber $(c_i,c_{i+1})$, the K-moduli stack $\MM_4^\K(c)$ is independent of $c$. For each wall $c_i$, let
\[
\phi_{c_i}^\pm\ \colon\ 
\ove{\fM}^{\K}(c_i\pm\varepsilon)
\ \longrightarrow\ 
\ove{\fM}^{\K}(c_i)
\]
denote the wall-crossing morphisms. We describe wall-crossings locally over their centers: a single wall may involve several distinct local wall-crossings, whose centers in our setting are isolated, pairwise disjoint points. For each such center $E_i\subseteq\ove{\fM}^{\K}(c_i)$, let $E_i^\pm$ denote the corresponding exceptional loci of $\phi_{c_i}^\pm$. When several local wall-crossings occur at the same wall, we distinguish their centers and exceptional loci using primed notation. Whenever $\phi_{c_i}^{\pm}$ are birational, write
\[
\psi_{c_i}\ \coloneqq \ (\phi_{c_i}^+)^{-1}\circ\phi_{c_i}^-
\ \colon\ 
\ove{\fM}^{K}(c_i-\varepsilon)
\ \dashrightarrow\ 
\ove{\fM}^{K}(c_i+\varepsilon)
\]
for the induced birational map. We say that the wall-crossing at $c_i$ is a \textit{divisorial contraction} if $\phi_{c_i}^-$ is an isomorphism and $\phi_{c_i}^+$ is a divisorial contraction, and a \textit{flip} if $\psi_{c_i}$ is a flip.

\subsection{K-moduli spaces of curves on $\bP^1\times \bP^1$} \label{K-moduli space of curves on smooth quadric surface}

We now restrict to K-moduli spaces that generically parametrize pairs $(\bP^1 \times \bP^1, cC)$, where $C$ is a smooth $(3,3)$-curve on $\bP^1\times \bP^1$.

\begin{prop}\label{interpolation}
Let $C$ be a $(3,3)$-curve on $\bP^1\times\bP^1$. If
$\operatorname{lct}(\bP^1\times\bP^1;C)>\tfrac23$
\textup{(}resp.\ $\geq\tfrac23$\textup{)}, then
$(\bP^1\times\bP^1,cC)$ is K-stable
\textup{(}resp.\ K-semistable\textup{)} for every $0<c<\tfrac23$.
In particular, it is K-stable if $C$ is smooth or has only
$A_n$-singularities with $n\leq4$, and K-semistable if $C$
has only $A_n$-singularities with $n\leq5$ and $D_4$-singularities.
\end{prop}

\begin{proof}
The first assertion follows from \Cref{Interpolation of K-stability ADL},
since $\bP^1\times\bP^1$ is K-polystable and
$(\bP^1\times\bP^1,\tfrac23C)$ is klt
\textup{(}resp.\ lc\textup{)}. The remaining assertions follow from $\lct(A_n)=\frac12+\frac1{n+1}$ and $\lct(D_n)=\frac{n}{2(n-1)}$.
\end{proof}

\begin{defn} \label{def: k moduli spaces of bidegree 3}
    For $c \in (0,\tfrac{2}{3}) \cap \bQ$, let $\MM_4^\K(c)$ denote the K-moduli stack of log Fano pairs $(X,cD)$ admitting a $\bQ$-Gorenstein smoothing to $(\bP^1\times \bP^1,cC)$, where $D \sim_{\bQ} -\tfrac{3}{2}K_{X}$. Let $\ove{\fM}_4^{K}(c)$ denote the good moduli space of $\MM_4^\K(c)$. 
\end{defn}

From Proposition \ref{interpolation}, we see that the spaces $\ove{\fM}_4^\K(c)$ are non-empty and birational to one another. This also justifies the subscript 4: a general smooth curve of genus 4 is realized as a smooth $(3,3)$-curve on $\bP^1 \times \bP^1$, so these spaces are birational to the moduli space $\overline{M}_4$ of DM stable curves of genus 4.

\subsection{VGIT for $(2,3)$-complete intersections in $\bP^{3}$}

K-moduli spaces in the small-coefficient regime are closely related to GIT and VGIT moduli spaces.

Let $\textbf{P}_{(3,3)}:= \bP(H^{0}(\bP^{1} \times \bP^{1}, \cO(3,3)))$. We say that a $(3,3)$-curve $C$ on $\bP^{1} \times \bP^{1}$ is \emph{GIT-(semi/poly)stable} if the point $[C]$ is GIT-(semi/poly)stable with respect to the natural $\Aut(\bP^{1} \times \bP^{1})$-action on $(\textbf{P}_{(3,3)}, \cO(2))$. We define the GIT moduli stack $\MM^{\GIT}$ and the GIT moduli space $\ove{\fM}^{\GIT}$ of $(3,3)$-curves by
\begin{equation*}
    \MM^{\GIT}\ \coloneqq \  [\textbf{P}_{(3,3)}^{\sst}/\Aut(\bP^{1} \times \bP^{1})], \qquad \ove{\fM}^{\GIT}\ \coloneqq \  \textbf{P}_{(3,3)}^{\sst}\sslash \big((\SL_2 \times \SL_2)\rtimes \bZ_2\big).
\end{equation*}

We recall the VGIT construction for $(2,3)$-complete intersections in $\bP^3$ from \cite{CMJL14}. Let
\[
\pi:\bP(E_3)\longrightarrow
\bP(H^0(\bP^3,\cO(2)))=\bP^9
\]
be the projective bundle whose fiber over a quadric $[Q]$ is
$\bP(H^0(Q,\cO_Q(3)))$. Write
\[
\eta=\pi^*\cO_{\bP^9}(1),
\qquad
\xi=\cO_{\bP(E_3)}(1),
\qquad
L_t=\eta+t\xi.
\]
The natural $\mathrm{SL}(4)$-action gives the slope-$t$ quotients
\[
\MM^{\GIT}(t)\coloneqq
[\bP(E_3)^{ss}(L_t)/\mathrm{PGL}(4)],
\qquad
\ove{\fM}^{\GIT}(t)\coloneqq
\bP(E_3)\sslash_{L_t}\mathrm{SL}(4),
\]
for $t\in(0,\tfrac23]\cap\bQ$.
Although $L_t$ is ample only for $0<t<\tfrac12$,
these quotients are well-defined for all such $t$,
and numerical \textup{(}semi\textup{)}stability agrees
with GIT \textup{(}semi\textup{)}stability by
\cite[Corollary 6.6]{CMJL14}.

The GIT moduli stack and space of $(3,3)$-curves coincide with the VGIT quotients in the first chamber:
\[
    \MM^{\GIT}(t)\ \simeq\ \MM^{\GIT},
    \qquad
    \ove{\fM}^{\GIT}(t)\ \simeq\ \ove{\fM}^{\GIT}
    \qquad\text{for }0<t<\tfrac29;
\]
see \cite{CMJL14}. Their relationship with K-moduli in the small-coefficient regime is described in \cite[Theorem 1.4]{ADL21}.

\subsection{Moduli of boundary polarized Calabi--Yau pairs}

In this section, we review the moduli theory of boundary polarized Calabi--Yau pairs developed in \cite{ABB+23, blum24}, in the setting of genus four curves. 

\begin{defn}
A \emph{boundary polarized CY pair} is a projective slc pair $(X,cD)$ such that 
\begin{enumerate}
    \item $D$ is an ample $\bQ$-Cartier Weil divisor on $X$ and $c\in \bQ_{>0}$;
    \item $K_X+cD\sim_{\bQ} 0$. 
\end{enumerate}
\end{defn}

\begin{theorem}[{\cite[Theorem 3.7]{ABB+23}}]
Let $N\in \bZ_{>0}$ and $c\in \bQ_{>0}$ with $Nc\in \bZ$, and let $\chi:\bN\to \bZ$ be a function. The moduli pseudo-functor sending each reduced scheme $B$ to
\[
\left\{\begin{array}{l}\textrm{families of boundary polarized CY pairs }(X,cD)\to B \textrm{ with coefficient $c$ and }\\ \textrm{index dividing $N$ satisfying }\chi(X_b, \omega_{X_b}^{[-m]}) = \chi(m) \textrm{ for all $b\in B$ and $m\in \bN$}
\end{array}\right\}
\] is represented by an algebraic stack $\MM^{\CY}(\chi, N, c)$, locally of finite type with affine diagonal.
\end{theorem}
In this paper we will focus on degenerations of pairs $(\bP^1\times\bP^1, \tfrac{2}{3}D)$ with $D\in |\cO_{\bP^1\times\bP^1}(3,3)|$.

\begin{defn}\label{def:CY4stack}
    Let $\chi(m) = \chi(\bP^1\times\bP^1, \omega_{\bP^1\times\bP^1}^{[-m]}) = (2m+1)^2$.
    Let $\MM_4^{\circ}$ denote the open substack of $\MM(\chi, 3,\tfrac{2}{3})$ whose $\bC$-points are the pairs $(\bP^1\times\bP^1, \tfrac{2}{3}D)$ with $D$ a smooth curve of bidegree $(3,3)$, and let $\mathcal{M}_4^{\CY}$ denote the stack-theoretic closure of $\mathcal{M}_4^{\circ}$ in $\MM(\chi, 3,\tfrac{2}{3})$. We call $\mathcal{M}_4^{\CY}$ the \emph{boundary polarized Calabi--Yau moduli stack} for genus four curves.
\end{defn}
\begin{theorem}[{\cite[Theorem 16.11]{ABB+23} and \cite[Theorem 7.5]{blum24}}]\label{thm:CY-gms}
    The algebraic stack $\mathcal{M}_4^{\CY}$ is irreducible and of finite type with affine diagonal, and it admits a good moduli space $\ove{\fM}_4^{\CY}$ as an irreducible projective scheme.
\end{theorem}

\begin{defn}\label{def:KAndKSBA}
Let $\mathcal{M}_4^{\K}$ (resp. $\mathcal{M}_4^{{\rm KSBA}}$) denote the substack of $\mathcal{M}_4^{\CY}$ parametrizing pairs $(X,\tfrac{2}{3}D)$ such that $(X, (\tfrac{2}{3}-\epsilon)D)$ is a K-semistable log Fano pair (resp. $(X, (\tfrac{2}{3}+\epsilon)D)$ is slc) for $0<\epsilon\ll 1$. 
\end{defn}

The $\mathcal{M}_4^{\K}$ agrees with $\mathcal{M}_4^{\K}( \frac{2}{3} - \epsilon)$ as defined in Theorem \ref{thm:Kmodthm} (see \cite[Remark 11.8]{ABB+23}).

\begin{theorem}[\cite{ABB+23, blum24}]\label{thm:bpCYwallcrossing}
    Both $\mathcal{M}_4^{\K}$ and $\mathcal{M}_4^{{\rm KSBA}}$ are open substacks of $\mathcal{M}_4^{\CY}$. Furthermore, there is a commutative wall-crossing diagram
        \[
\begin{tikzcd}
\mathcal{M}_4^{\K} \arrow[d] \arrow[r,hook]  &
\mathcal{M}_4^{\CY} \arrow[d] & 
\mathcal{M}_4^{{\rm KSBA}} \arrow[d] \arrow[l,hook']  \\
\ove{\fM}_4^{\K} \arrow[r]  & 
\ove{\fM}_4^{\CY}  & 
\ove{\fM}_4^{{\rm KSBA}} \arrow[l]
\end{tikzcd}
,\]
in which the top arrows are open immersions, the vertical arrows are good moduli space morphisms, and the bottom arrows are projective birational morphisms. 
\end{theorem}

\begin{proof}
This follows from \cite[Theorems 1.1 and 7.5]{blum24} (see also \cite[Theorem 16.11]{ABB+23}). The normality of $\ove{\fM}_4^{\K}$ follows from \cite[Theorem 2.21]{ADL21}.
\end{proof}

Finally, we observe that these moduli spaces are all birational; the following result shows that the locus of klt pairs $(X,\tfrac{2}{3}D)$ is contained in all three stacks.  

\begin{prop}\cite[Proposition 2.34]{ADLW25}\label{prop:klt-is-stable}
    If $[(X,\tfrac{2}{3}D)] \in \mathcal{M}_4^{\CY}$ is a klt pair, then $[(X,\tfrac{2}{3}D)] \in  \mathcal{M}_4^{{\KSBA}}$ and $[(X,\tfrac{2}{3}D)] \in \mathcal{M}_4^{\K}$. In fact, $(X,(\tfrac{2}{3}- \epsilon)D)$ is K-stable for all $0<\epsilon \ll 1$.
\end{prop}

\subsection{Classification of degenerations of $\bP^{1} \times \bP^{1}$}
In this final section, we describe all possible degenerations of $\bP^{1} \times \bP^{1}$ appearing in $\MM^\K_4(c)$.  We first introduce the surface $X_{10}$, a partial smoothing of $\bP(1,2,9)$.

\begin{defn}\label{def:X10}
  Let $[x:y:z:w]$ be the weighted homogeneous coordinates of $\bP(1,1,5,9)$, and define the weighted hypersurface
\[
X_{10} \ \coloneqq\ 
V(xw-y^{10}-z^2)\ \subseteq \ \bP(1,1,5,9).
\]  
This has a unique singular point at $[0:0:0:1]$ of type $\frac{1}{9}(1,2)$.  Observe that $\bP(1,2,9) \cong V(xw - z^{2})$ and hence $X_{10}$ is a partial smoothing of $\bP(1,2,9)$.
\end{defn}

\begin{prop}\label{new surfaces in K-moduli spaces}
Let $(X,cD)$ be a K-semistable log Fano pair admitting a
$\bQ$-Gorenstein smoothing to $(\bP^1\times\bP^1,cC)$,
where $0<c<\tfrac23$ and $C$ is a $(3,3)$-curve.
Then $X$ is isomorphic to $\bP^1\times\bP^1$,
$\bP(1,1,2)$, $X_{10}$, or $\bP(1,2,9)$.
\end{prop}

\begin{proof}
If $X$ is smooth, then $X\simeq\bP^1\times\bP^1$.
If $\rho(X)=1$, the argument of
\cite[Theorem 4.8]{ADL21} gives
$X\simeq\bP(1,1,2)$ or $\bP(1,2,9)$.

Suppose that $\rho(X)=2$ and $X$ is singular.
By \cite[Proposition 2.6]{HP10}, its singularities
are of type $\frac1{n^2}(1,na-1)$ with $2\leq n\leq4$.
The case $n=4$, corresponding to
$\frac1{16}(1,3)$, is excluded by the argument of
\cite[Theorem 4.8]{ADL21}. Thus only
$\frac14(1,1)$ and $\frac19(1,2)$ can occur. By \cite[Proposition 3.1]{HP10}, the $\bQ$-Gorenstein deformations of $X$ are unobstructed and allow its singularities to be smoothed independently.  Furthermore, the K-moduli stack $\MM^\K_4(c)$ is smooth by the argument of \cite[Proposition 11.7]{ABB+23}, and hence by openness of K-semistability in families \cite{BLX19}, the generic smoothing will remain a K-semistable pair that admits a smoothing to $\bP^1 \times \bP^1$. If $X$ has a $\frac14(1,1)$-singularity, smooth all other singularities. By \cite[Remark 4.4]{Fujita2016LogDelPezzo}, the resulting
surface is the blowup of $\bP(1,1,4)$ at a smooth
point. Its $\bQ$-Gorenstein smoothing is $\bF_1$, a contradiction. Hence $X$ has only $\frac19(1,2)$-singularities, and \cite[Remark 4.4]{Fujita2016LogDelPezzo}
gives $X\simeq X_{10}$.
\end{proof}

\section{GIT of $(3,3)$-curves on $\bP^1\times\bP^1$ and the first wall} \label{sec: first wall}

By \cite[Theorem 1.4]{ADL21}, for $c \in (0, \frac{1}{6})$ we have an isomorphism $\ove{\fM}_4^\K(c) \simeq \ove{\fM}^{\GIT}$ where the latter is the GIT moduli space of $(3,3)$-curves on $\bP^1 \times \bP^1$. In this section, we collect some results on the GIT classification, and study the geometry of the first K-moduli wall $c_{1} = \tfrac{1}{6}$ and the relation between the K-moduli space and the GIT moduli space of $(3,3)$-curves on $\bP^1\times\bP^1$.

\subsection{GIT of $(3,3)$-curves on $\bP^{1} \times \bP^{1}$} \label{section 4}

We begin by reviewing the GIT quotient of $V:= |\mathcal{O}_{\bP^{1} \times \bP^{1}}(3,3)|$ by $\Aut(\bP^1\times\bP^1)$. For more details, we refer the reader to \cite{Fed12, CMJL14}.

\begin{prop}[cf.\ {\cite[Prop.~2.2]{Fed12}}] \label{unstable GIT (3,3) curve}
    A $(3,3)$-curve $C$ on $\bP^1\times \bP^1$ is GIT unstable if and only if one
of the following holds:
\begin{enumerate}
    \item $C$ contains a double ruling; or
    \item $C$ contains a ruling $R$ that meets the residual curve $C' = \overline{C - R}$ at a single point, which is also a singular point of $C'$.
\end{enumerate}
\end{prop}

The intersection of the smooth quadric with a hyperplane in $\bP^3$, when smooth, is a plane conic; accordingly, we refer to smooth $(1,1)$-curves on $\bP^1 \times \bP^1$ as \textit{conics}.

\begin{remark}\label{non-reduced GIT semistable curve}
By Proposition \ref{unstable GIT (3,3) curve} (cf.\ \cite[Corollary 2.3]{Fed12}), if $C$ is a non-reduced GIT semistable $(3,3)$-curve, then $C$ is one of the following:
\begin{enumerate}
    \item a triple conic; or
    \item the union of a smooth double conic and a reduced $(1,1)$-curve meeting the double conic transversely; or
    \item the union of a smooth double conic and a reduced $(1,1)$-curve meeting the double conic at a single point with multiplicity $2$; in this case the reduced $(1,1)$-curve is necessarily smooth.
\end{enumerate}
The curves of types (1) and (2) are GIT polystable but not stable. The curves of type (3) are strictly semistable and contain the triple conic in their orbit closure. 
\end{remark}

Next, we classify the reduced GIT semistable $(3,3)$-curves. By \cite[Proposition 4.11]{CMJL14}, if a reduced $(3,3)$-curve has a singularity of multiplicity greater than three, then it is GIT unstable.

\begin{prop}\textup{(cf.\ \cite[Prop.~2.1, 2.5]{Fed12})} \label{prop:D-curve}
    Let $C$ be a reduced GIT semistable $(3,3)$-curve with a singularity $p$ of multiplicity $3$. Then the orbit closure of $C$ contains the point representing the tangent cone of $C$ at $p$, which is a union of three conics meeting in two $D_{4}$-singularities.
\end{prop}

Recall that a planar curve singularity of multiplicity two is an $A_{k}$-singularity, locally analytically of the form $x^2 + y^{k+1} = 0$. We say that an $A_k$-singularity of a curve $C$ is \textbf{separating} if the partial normalization $\widetilde{C} \to C$ resolving only this singularity is disconnected, and \textbf{non-separating} if $\widetilde{C}$ is connected. This distinction is relevant only for odd $k$.

\begin{prop}\textup{(cf.\ \cite[Prop.~2.1, 3.6]{CMJL14})} \label{prop:Ak singularity of (3,3) curve}
    Let $C$ be a reduced $(3,3)$-curve on $\bP^1\times\bP^1$ with an $A_k$-singularity $p$. If $p$ is non-separating, then $k \leq 8$. If $p$ is separating, then one of the following holds:
\begin{enumerate}
    \item $k = 9$, and $C$ is a union of two twisted cubics; or
    \item $k = 7$, and $C$ is a union of a quartic and a conic; or
    \item $k = 5$, and $C$ is a union of a quintic and a line.
\end{enumerate}
Moreover, every $(3,3)$-curve with only double points is GIT semistable.
\end{prop}

By Proposition \ref{interpolation} and \cite[Theorem 1.4]{ADL21}(1), \cite[Proposition 4.7]{ADL24} provides open embeddings
\begin{equation*}
    \fM^{\GIT,\klt} \ \hookrightarrow \ \fM^{\GIT,\lc} \ \hookrightarrow \  \ove{\fM}^{\GIT} \cong \ove{\fM}_4^\K(\tfrac{1}{6} - \epsilon), 
\end{equation*}
where $\fM^{\GIT,\klt}$ and $\fM^{\GIT,\lc}$ denote the loci in the GIT moduli space $\ove{\fM}^{\GIT}$ parametrizing the GIT-polystable $(3,3)$-curves with $\operatorname{lct} > \tfrac{2}{3}$ and $\operatorname{lct} \geq \tfrac{2}{3}$, respectively.

\begin{lemma}\label{strata-for-GIT}
    The boundary $\ove{\fM}^{\GIT} \setminus \fM^{\GIT,\klt}$ is a disjoint union of the following locally closed strata.

    \emph{Zero-dimensional loci:}
    \begin{itemize}
\item $\Sigma_1$, parametrizing a triple conic;
\item $\Sigma_2$, parametrizing a $(3,3)$-curve with a separating $A_{9}$-singularity; and
\item $\Sigma_7$, parametrizing a maximally degenerate $(3,3)$-curve with two $A_{5}$-singularities, consisting of two lines in the same ruling $(1,0)$ and a smooth $(1,3)$-curve meeting each line at an $A_5$-singularity (i.e., tangent to order $3$).
\end{itemize}

    \emph{One-dimensional loci:}
    \begin{itemize}
\item $\Sigma_3$, parametrizing $(3,3)$-curves with an $A_{8}$-singularity; and
\item $\Sigma_9$, parametrizing unions of a double conic and another conic.
\end{itemize}

    \emph{Two-dimensional loci:}
    \begin{itemize}
\item $\Sigma_4$, parametrizing $(3,3)$-curves with a separating $A_{7}$-singularity;
\item $\Sigma_5$, parametrizing $(3,3)$-curves with a non-separating $A_{7}$-singularity; and
\item $\Sigma_{10}$, parametrizing $D$-curves, which are unions of three conics meeting in two $D_4$ singularities
\end{itemize}

    \emph{Three-dimensional locus:}
    \begin{itemize}
\item $\Sigma_6$, parametrizing $(3,3)$-curves with an $A_{6}$-singularity.
\end{itemize}

    \emph{Four-dimensional locus:}
    \begin{itemize}
    \item $\Sigma_8$, parametrizing $(3,3)$-curves with a non-separating $A_5$-singularity.
\end{itemize}
Moreover, we have the following incidence diagram:

\begin{equation*}
\begin{tikzcd}
    \Sigma_7 & & \arrow[ld] \Sigma_4 &   &  \arrow[ld] \Sigma_2 \\
    \Sigma_8 & \arrow[l] \Sigma_6 & \arrow[l] \Sigma_5 & \arrow[l] \Sigma_3 & \arrow[l] \arrow[ld] \Sigma_1 \\
     & & \Sigma_{10} & \arrow[l] \Sigma_9 &
\end{tikzcd}
\end{equation*}
Here $\Sigma_{j} \leftarrow \Sigma_{i}$ means that $\Sigma_{i}$ is contained in the boundary of the Zariski closure of $\Sigma_{j}$ in $\ove{\fM}^{\GIT}$. Furthermore, $\Sigma_{2} \not\subseteq \overline{\Sigma}_4$, and the boundary of $\Sigma_4$ meets $\overline{\Sigma}_9$ in a single point, representing the triple conic. The locus $\Sigma_1 \cup \Sigma_7 \cup \Sigma_9 \cup \Sigma_{10}$ is exactly the strictly GIT polystable locus, and the remaining strata lie in the stable locus. In addition, $\ove{\fM}^{\GIT} \setminus \fM^{\GIT,\lc} = \bigcup_{i=1}^{6}\Sigma_{i} \cup \Sigma_{9}$. 
\end{lemma}

\begin{proof}
    By the Hilbert--Mumford numerical criterion \cite{MFK94} and the GIT analysis in \cite{Fed12}, a reduced $(3,3)$-curve $C$ satisfies the following.
    \begin{itemize}
        \item $C$ is strictly GIT polystable if and only if $C$ is either a $D$-curve or a maximally degenerate $(3,3)$-curve with two $A_{5}$-singularities.
        \item $C$ is GIT stable if and only if $C$ is smooth or all of its singularities are of type $A_{k}$ with $1 \leq k \leq 9$ and $C$ has no separating $A_5$-singularity.
    \end{itemize}
    Since the GIT moduli space $\ove{\fM}^{\GIT} $ is the union of the GIT stable locus and the strictly GIT polystable locus, the lemma follows by analyzing the log canonical thresholds of these singularities, together with the relevant jet computations for the dimensions of the strata (see Sections \ref{sec:walls of P(1,1,2)}, \ref{sec:walls of X10}, and \ref{sec: walls of P(1,2,9)}).
    
    We observe that the general member of $\Sigma_4$ has a component of degree $2$, whereas a curve in $\Sigma_2$ is a union of two twisted cubics and thus has no such component; hence $\Sigma_2$ is not contained in the boundary of $\Sigma_4$.

    We also consider points in the intersection of $\overline{\Sigma_9}$ and $\overline{\Sigma_4}$. Such a point must be a union of a double conic and another conic having an intersection point of order at least four. The only GIT semistable curves with this property are those of types (1) and (3) in Remark \ref{non-reduced GIT semistable curve}, and therefore the only polystable curve with this property is the triple conic.
\end{proof}

\begin{remark} \label{singularities in wall crossings}
   By \cite[Theorem 1.4]{ADL21} and \Cref{interpolation}, every curve $C$ parametrized by $\ove{\fM}^{\GIT,\klt}$ gives a pair $(\bP^1 \times \bP^1, cC)$ that is K-stable for all $c \in (0, \tfrac{2}{3})$. Therefore, to compute the wall crossings, it suffices to consider the non-klt curves, i.e. the $(3,3)$-curves with the following singularities: triple conics, double conics, non-separating $A_{k}$-singularities with $6\leq k \leq 8$, and separating $A_{k}$-singularities with $k \in \{7, 9\}$.
\end{remark}

\subsection{The first wall crossing}
\label{subsec:relation to GIT}

In this subsection we identify the K-moduli space for
small coefficient with the GIT moduli space of
$(3,3)$-curves and describe the first wall crossing.

Set $c_1=\tfrac16$. Let $H\subset\bP^1\times\bP^1$
be a smooth $(1,1)$-curve, and set
\[
C_{\mathrm{tc}}\coloneqq3H,\qquad
C'_{\mathrm{tc}}\coloneqq V(s^3)\subset\bP(1,1,2).
\]

\begin{theorem}\label{thm:first-wall}
Keep the notation introduced above. Then
\begin{enumerate}
\item[\textup{(1)}]
For every $0<c<c_1$, a pair
$(\bP^1\times\bP^1,cC)$ is K-(poly/semi)stable
if and only if $[C]$ is GIT-(poly/semi)stable.
Moreover, there is an isomorphism
\[
\ove{\fM}^{\K}_4(c)\simeq\ove{\fM}^{\GIT}.
\]

In general, if $(\bP^1 \times \bP^1, cC)$ is K-(poly/semi)stable for some $c \in (0, \frac{2}{3})$, then $C$ is GIT-(poly/semi)stable. 

\item[\textup{(2)}]
The pair $(\bP^1\times\bP^1,c_1C_{\mathrm{tc}})$
is K-semistable, with K-polystable degeneration
$(\bP(1,1,2),c_1C'_{\mathrm{tc}})$.
Up to automorphisms of $\bP(1,1,2)$,
$C'_{\mathrm{tc}}$ is the unique curve of weighted
degree $6$ such that the pair $(\bP(1,1,2), c_1 C)$
is K-polystable.

\item[\textup{(3)}]
The wall-crossing morphism
\[
\phi^-_{c_1}\colon
\ove{\fM}^{\K}_4(c_1-\epsilon)
\longrightarrow\ove{\fM}^{\K}_4(c_1)
\]
is an isomorphism and identifies the triple conic point
with $[(\bP(1,1,2),c_1C'_{\mathrm{tc}})]$.
The morphism
\[
\phi^+_{c_1}\colon
\ove{\fM}^{\K}_4(c_1+\epsilon)
\longrightarrow\ove{\fM}^{\K}_4(c_1)
\]
is a weighted blowup at this point.
Its exceptional divisor parametrizes GIT-polystable
curves of weighted degree $6$ on $\bP(1,1,2)$
in the sense of \Cref{defn:gitp112-1}
and is birational to the Gieseker--Petri divisor.\footnote{Observe that $\ove{\fM}^{\K}_4(c)$ is birational to $\overline{M}_4$, the moduli space of stable curves of genus 4.  Via this birational map, the exceptional divisor of the morphism $\phi^+_{c_1}$ is the birational transform of the Gieseker--Petri divisor on $\overline{M}_4$.} 
\end{enumerate}
\end{theorem}

\begin{proof}[Proof of \Cref{thm:first-wall}\textup{(1)}]
The proof is the same as that of
\cite[Theorem 5.2]{ADL21}.
\end{proof}

To prove the remaining assertions, we first exclude
curves passing through the vertex of $\bP(1,1,2)$.

\begin{lemma}\label{lem:no-nodes-on-cone-before-1/3}
Let $D\subset\bP(1,1,2)$ be a curve of weighted
degree $6$. If $D$ passes through the vertex, then
$(\bP(1,1,2),cD)$ is K-unstable for every
$0<c<\tfrac13$.
\end{lemma}

\begin{proof}
Let $\mu:Y\to\bP(1,1,2)$ be the minimal resolution,
and let $E$ be its exceptional divisor.
On the index cover of the vertex, the pullback of
$D$ has multiplicity at least $2$. Thus
$\ord_E(D)\geq1$, and
\[
A_{(\bP(1,1,2),cD)}(E)\leq1-c.
\]
A direct computation gives
\[
S_{(\bP(1,1,2),cD)}(E)
=\tfrac{4}{3}\left(1-\tfrac{3}{2}c\right).
\]
Consequently,
\[
\beta_{(\bP(1,1,2),cD)}(E)
\leq-\tfrac{1}{3}+c<0
\qquad\text{for }0<c<\tfrac{1}{3}.
\]
The assertion follows from the valuative criterion.
\end{proof}

\begin{proof}[Proof of \Cref{thm:first-wall}\textup{(2)}]
The construction of the degeneration and the
K-polystability of
$(\bP(1,1,2),c_1C'_{\mathrm{tc}})$ follow by the
same argument as in \cite[Proposition 5.3]{ADL21}.
Openness of K-semistability then implies that
$(\bP^1\times\bP^1,c_1C_{\mathrm{tc}})$
is K-semistable.

For the uniqueness assertion, let
$(\bP(1,1,2),c_1D)$ be K-polystable.
By \Cref{lem:no-nodes-on-cone-before-1/3},
the curve $D$ avoids the vertex. After a change
of coordinates, its equation can be written as
\[
s^3+f_4(u,v)s+f_6(u,v)=0,
\]
where each $f_i$ is a binary form of degree $i$.
The family
\[
s^3+t^2f_4(u,v)s+t^3f_6(u,v)=0
\]
gives a $\bG_m$-equivariant degeneration to
$C'_{\mathrm{tc}}$.
The central pair is K-polystable, so this is a
degeneration between K-semistable pairs with
vanishing Futaki invariant.
K-polystability therefore implies
\[
(\bP(1,1,2),D)\simeq
(\bP(1,1,2),C'_{\mathrm{tc}}).
\]
\end{proof}

We next define the GIT stability condition that will define
the curves on $\bP(1,1,2)$ appearing after the first wall.

\begin{defn}\label{defn:gitp112-1}
Let $D\subset\bP(1,1,2)$ be a curve of weighted
degree $6$ avoiding the vertex.
After a change of coordinates, its equation
can be written as
\begin{equation}\label{eq:even-p112}
s^3+f_4(u,v)s+f_6(u,v)=0.
\end{equation}
Set
\[
\bA'_6\coloneqq
H^0(\bP^1,\cO(4))\oplus H^0(\bP^1,\cO(6)),
\]
and let $\bG_m$ act by
\[
\lambda\cdot(f_4,f_6)
=(\lambda^2f_4,\lambda^3f_6).
\]
Let $\bP'_6$ be the coarse moduli space of
$[(\bA'_6\setminus\{0\})/\bG_m]$.

The natural $\SL(2)$-action on binary forms
commutes with this $\bG_m$-action and descends
to $\bP'_6$.
For $(f_4,f_6)\neq(0,0)$, we say that
$[D]\in\bP'_6$ is \emph{GIT-(poly/semi)stable}
if it is GIT-(poly/semi)stable for this
$\SL(2)$-action with respect to
$\cO_{\bP'_6}(1)$, viewed as a $\bQ$-line bundle.
\end{defn}

\begin{theorem}\label{thm:k=gitp112-1}
Let $D\subset\bP(1,1,2)$ be a curve of weighted
degree $6$ represented by a point $[D]\in\bP'_6$.  If $(\bP(1,1,2), cD)$ is K-(poly/semi)stable for some $c \in (c_1, \frac{2}{3})$, then $D$ is GIT-(poly/semi)stable in the sense of \Cref{defn:gitp112-1}.

For $0<\epsilon\ll1$, the pair
$(\bP(1,1,2),(c_1+\epsilon)D)$ is
K-(poly/semi)stable if and only if $[D]$
is GIT-(poly/semi)stable in the sense of
\Cref{defn:gitp112-1}.  Furthermore, all such pairs admit an isotrivial specialization to $(\bP(1,1,2), C_{\mathrm{tc}}')$.
\end{theorem}

\begin{proof}
The proof follows the argument of
\cite[Theorem 5.8]{ADL24}.
\end{proof}

\begin{proof}[Proof of \Cref{thm:first-wall}\textup{(3)}]
The proof follows from \cite[Theorem 1.4]{ADL21} and the description of the VGIT wall-crossings in \cite{CMJL14}, although we provide an alternative proof below that does not rely on comparison with VGIT. 

The wall-crossing morphisms fit into a diagram
\[
\begin{tikzcd}
\MM^{\K}(c_1-\epsilon)
  \arrow[r]\arrow[d]
&
\MM^{\K}(c_1)
  \arrow[d]
&
\MM^{\K}(c_1+\epsilon)
  \arrow[l]\arrow[d]
\\
\ove{\fM}^{\GIT}
  \simeq\ove{\fM}^{\K}_4(c_1-\epsilon)
  \arrow[r,"\phi^-_{c_1}"]
&
\ove{\fM}^{\K}_4(c_1)
&
\ove{\fM}^{\K}_4(c_1+\epsilon).
  \arrow[l,"\phi^+_{c_1}"']
\end{tikzcd}
\]

We claim that $\phi^-_{c_1}$ must be an isomorphism.  First note that the K-moduli spaces are normal by
\cite[Theorem 2.21]{ADL21}, so $\phi^-_{c_1}$ is a birational morphism of normal varieties.  Second, by the same argument as in \cite[Theorem 5.5]{ADL24}, the Picard rank of $\ove{\fM}^{\GIT}$ is 1, so we conclude $\phi^-_{c_1}$ is finite.  Therefore, $\phi^-_{c_1}$ is an isomorphism.  

From this, we claim the polystable points of $\ove{\fM}^{\K}_4(c_1)$ are precisely the pairs $(\bP^1 \times \bP^1, c_1 C)$ such that $C$ is GIT polystable and not projectively equivalent to $C_{\mathrm{tc}}$ and the pair $(\bP(1,1,2), c_1 C'_{\mathrm{tc}})$.  Suppose $C$ is a GIT polystable curve and $(X,c_1D)$ is the K-polystable degeneration of $(\bP^1 \times \bP^1, c_1C)$.  By \cite[Theorem 4.7]{ADL21}, we conclude $X = \bP^1 \times \bP^1$ or $X = \bP(1,1,2)$.  If $X = \bP(1,1,2)$, then by Theorem \ref{thm:first-wall}(2), we must have $(X,c_1D) = (\bP(1,1,2), c_1 C'_{\mathrm{tc}})$, i.e. $\phi^-_{c_1}([C]) = \phi^-_{c_1}([C_{\mathrm{tc}}])$, contradicting the injectivity of $\phi^-_{c_1}$.  Therefore, $X = \bP^1 \times \bP^1$ and the GIT polystability of $C$ implies $(X,c_1D) \cong (\bP^1 \times \bP^1, c_1C)$ so $(\bP^1 \times \bP^1, c_1 C)$ is K-polystable.  The result now follows from Theorem \ref{thm:first-wall}(2) and this completes the description of $\phi^-_{c_1}$.  

To describe $\phi^+_{c_1}$, we first claim that the K-polystable pairs $(X, (c_1+\epsilon)D)$ are the pairs $(\bP^1 \times \bP^1, (c_1+\epsilon) C)$ such that $C$ is GIT polystable and not projectively equivalent to $C_{\mathrm{tc}}$ and pairs $(\bP(1,1,2), (c_1+\epsilon)D)$ such that $D$ is GIT polystable in the sense of Definition \ref{defn:gitp112-1}.   As above, by \cite[Theorem 4.7]{ADL21}, we must have $X = \bP^1 \times \bP^1$ or $X = \bP(1,1,2)$.  Suppose $X = \bP^1 \times \bP^1$.  If $(\bP^1 \times \bP^1, (c_1+\epsilon)C)$ is K-polystable, then $C$ is GIT polystable by the same proof as \cite[Theorem 1.4, Theorem 5.2]{ADL21}.  Furthermore, $(\bP^1 \times \bP^1, (c_1 + \epsilon) C_{\mathrm{tc}})$ is K-unstable for any $\epsilon > 0$; this follows from interpolation and the fact that $(\bP^1 \times \bP^1, c_1 C_{\mathrm{tc}})$ is not K-polystable.  Therefore, the K-polystability of the pair $(\bP^1 \times \bP^1, (c_1+\epsilon)C)$ implies that $C$ is GIT polystable and not projectively equivalent to the triple conic.  Conversely, if $C$ is GIT polystable and not projectively equivalent to $C_{\mathrm{tc}}$, then the previous paragraph demonstrates that $(\bP^1 \times \bP^1, cC)$ is K-polystable for all $c\in (0, c_1]$ so it follows from \cite[Proposition 3.18]{ADL24} that $(\bP^1 \times \bP^1, (c_1 + \epsilon)C)$ is K-polystable.  Therefore, the K-polystable pairs $(\bP^1 \times \bP^1, (c_1 + \epsilon)C)$ are precisely the pairs where $C$ is GIT polystable and not projectively equivalent to the triple conic.  The case when $X = \bP(1,1,2)$ follows directly from Theorem \ref{thm:k=gitp112-1}.

Finally, the local VGIT presentation, combined with the
argument of \cite[Section 5.2]{ADL24} and
\cite{DH98}, identifies $\phi^+_{c_1}$ as a
weighted blowup at
$[(\bP(1,1,2),c_1C'_{\mathrm{tc}})]$.
By \Cref{thm:k=gitp112-1}, its exceptional divisor
parametrizes GIT-polystable curves of weighted
degree $6$ on $\bP(1,1,2)$ in the sense of
\Cref{defn:gitp112-1}.
The smooth genus 4 curves on $\bP(1,1,2)$, which are GIT-polystable in this definition, are precisely the Gieseker--Petri special canonical curves so this is the birational transform of the
Gieseker--Petri divisor.
\end{proof}

\section{Wall crossings involving $\bP(1,1,2)$}
\label{sec:walls of P(1,1,2)}

In this section we study the wall crossings involving
pairs whose underlying surface is $\bP(1,1,2)$.
The corresponding VGIT walls are associated with
triple conics, double conics, and $(3,3)$-curves with
a separating $A_7$-singularity; see \cite{CMJL14,ADL21}.
Having treated triple conics in the previous section,
we now consider the remaining two types.

\subsection{Double conics}

Let $C_{\mathrm{dc}}=2H+Q$, where $H$ and $Q$ are
smooth $(1,1)$-curves meeting transversely.
By \cite{Fed12}, the double conic locus in
$\ove{\fM}^{\GIT}$ is one-dimensional.

Set $c_2=\tfrac13$, and let
\[
C'_{\mathrm{dc}}\coloneqq V(uvs^2)\subset\bP(1,1,2).
\]
Thus $C'_{\mathrm{dc}}=2H_0+Q_0$, where
$H_0=V(s)$ is the section at infinity and
$Q_0=V(uv)$ has a node at the vertex.

\begin{theorem}\label{Second wall crossing}
Keep the notation introduced above. Then
\begin{enumerate}
\item[\textup{(1)}]
The second K-moduli wall occurs at $c_2=\tfrac13$.
The pair $(\bP^1\times\bP^1,cC_{\mathrm{dc}})$
is K-semistable \textup{(}resp.\ K-polystable\textup{)}
if and only if $c\leq c_2$
\textup{(}resp.\ $c<c_2$\textup{)}.

\item[\textup{(2)}]
The K-polystable degeneration of
$(\bP^1\times\bP^1,c_2C_{\mathrm{dc}})$
is isomorphic to
$(\bP(1,1,2),c_2C'_{\mathrm{dc}})$.

\item[\textup{(3)}]
The birational map
\[
\phi_{c_2}\colon
\ove{\fM}^{\K}_4(c_2-\epsilon)
\dashrightarrow
\ove{\fM}^{\K}_4(c_2+\epsilon)
\]
is a flip, with the exceptional loci of
$\phi^+_{c_2}$ and $\phi^-_{c_2}$ having dimensions
$7$ and $1$, respectively.
The former parametrizes GIT-polystable curves
\[
D=V\bigl(uvs^2+f_4(u,v)s+f_6(u,v)\bigr)
\subseteq\bP(1,1,2)
\]
in the sense of \Cref{defn:gitp112-2}, and the latter is the strict transform of the locus $\Sigma_9$.
\end{enumerate}
\end{theorem}

\begin{proof}[Proof of
\Cref{Second wall crossing}\textup{(1)--(2)}]
We first show that
$(\bP^1\times\bP^1,cC_{\mathrm{dc}})$
is K-unstable for $c>c_2$. We have
\[
A_{(\bP^1\times\bP^1,cC_{\mathrm{dc}})}(H)
=1-2c,
\qquad
S_{(\bP^1\times\bP^1,cC_{\mathrm{dc}})}(H)
=\tfrac{2}{3}-c.
\]
Therefore,
\[
\beta_{(\bP^1\times\bP^1,cC_{\mathrm{dc}})}(H)
=\tfrac{1}{3}-c<0
\qquad\text{for }c>c_2.
\]

We will show the rest of (1) and (2) via the following lemma.

\begin{lem}\label{de-double conic}
The pair $(\bP^1\times\bP^1,cC_{\mathrm{dc}})$
admits a special degeneration to
$(\bP(1,1,2),cC'_{\mathrm{dc}})$.
\end{lem}

\begin{proof}
Let $\sX=(\bP^1\times\bP^1)\times\bA^1_t$,
and let $\pi:\sX'\to\sX$ be the blowup along
$H\times\{0\}$.
Its central fiber is $X\cup S$, where
$X\simeq\bP^1\times\bP^1$ is the strict transform
of the original central fiber and $S\simeq\bF_2$
is the exceptional divisor.
The intersection $F=X\cap S$ is the negative
section of $S$.

Consider the relative divisor
\[
L=-K_{\sX'}-S.
\]
Its restriction to a general fiber is in the linear system $|\cO(2,2)|$,
while
\[
L|_X\sim_{\bQ}0,
\qquad
L|_S\sim_{\bQ}-K_{\bF_2}.
\]
The relative ample model of $L$ contracts $X$ to a point and restricts on $S$ to the contraction of its negative section. The resulting central fiber is therefore $\bP(1,1,2)$.

The limit of $H$ is the section at infinity $H_0$. Since $Q$ meets $H$ transversely at two distinct points, its transform on $S$ consists of two distinct ruling fibers. Their images form $Q_0=V(uv)$, which has a node at the vertex. Thus the limiting boundary is
$2H_0+Q_0=C'_{\mathrm{dc}}$.
\end{proof}

The degeneration in \Cref{de-double conic} is induced by $\ord_H$. The pair $(\bP(1,1,2),c_2C'_{\mathrm{dc}})$
is K-polystable by \cite[Theorem 6.3]{CMJL14} and \cite[Theorem 1.4]{ADL21}.
Openness of K-semistability implies that
$(\bP^1\times\bP^1,c_2C_{\mathrm{dc}})$
is K-semistable. Together with \Cref{thm:first-wall}\textup{(1)}
and interpolation of K-stability, this gives K-polystability for $0<c<c_2$.
The degeneration at $c_2$ is nontrivial, so the original pair is not K-polystable at the wall.

It remains to exclude walls in
$(c_1,c_2)=(\tfrac16,\tfrac13)$.
The classification and comparison results of
\cite{CMJL14,ADL21} show that a new K-polystable
replacement in this interval would have underlying
surface $\bP(1,1,2)$ and boundary passing through
the vertex.
This is excluded by
\Cref{lem:no-nodes-on-cone-before-1/3}.
Hence $c_2$ is the second wall.
\end{proof}

To complete the description of the wall crossing in (3), we next describe the curves introduced after $c_2$.

\begin{prop}\label{curves meet the vertex point in a node}
Let $(\bP(1,1,2),(c_2+\epsilon)D)$ be K-semistable,
where $0<\epsilon\ll1$.
If $D$ passes through the vertex, then it has a
node there. A general such curve is smooth away
from the vertex.
\end{prop}

\begin{proof}
This follows from \cite[Proposition 4.12]{CMJL14}
and the K/VGIT comparison in \cite{ADL21}.
\end{proof}

\begin{prop}\label{curves meet the vertex point in a cusp}
Let $D\subset\bP(1,1,2)$ be a curve of weighted
degree $6$ passing through the vertex $p$.
Then
\begin{enumerate}
\item[\textup{(1)}]
if $p$ is a cusp of $D$, then
$(\bP(1,1,2),cD)$ is K-unstable for
$0<c<\tfrac12$;
\item[\textup{(2)}]
if $p$ is neither a node nor a cusp of $D$, then
$(\bP(1,1,2),cD)$ is K-unstable for
$0<c<\tfrac23$.
\end{enumerate}
\end{prop}

\begin{proof}
This follows from \cite[Proposition 4.12]{CMJL14}
and the K-moduli and VGIT quotient comparison in \cite{ADL21}.
\end{proof}

\begin{lem}\label{lem:no-other-wall112}
Let $D\subseteq\bP(1,1,2)$ be a curve of weighted degree $6$, with $f_4$ and $f_6$ binary forms of
degrees $4$ and $6$, respectively.
Then
\begin{enumerate}
\item[\textup{(1)}]
if $D=V(uvs^2+f_4s+f_6)$, then
$(\bP(1,1,2),cD)$ is K-unstable for
$0<c<\tfrac13$;
\item[\textup{(2)}]
if $D=V(u^2s^2+f_4s+f_6)$, then
$(\bP(1,1,2),cD)$ is K-unstable for
$0<c<\tfrac12$;
\item[\textup{(3)}]
if $D=V(f_4s+f_6)$, then
$(\bP(1,1,2),cD)$ is K-unstable for
$0<c<\tfrac23$.
\end{enumerate}
\end{lem}

\begin{proof}
Assertion \textup{(1)} follows from
\Cref{lem:no-nodes-on-cone-before-1/3}.
In \textup{(2)}, the curve is not nodal at the vertex of $\bP(1,1,2)$ so the assertion follows from
\Cref{curves meet the vertex point in a cusp}.
In \textup{(3)}, the singularity is neither a node nor a cusp at the vertex and the same
proposition gives the result.
\end{proof}

\begin{defn}\label{defn:gitp112-2}
Let $D\subset\bP(1,1,2)$ be a curve with equation
\[
uvs^2+f_4(u,v)s+f_6(u,v)=0.
\]
After replacing $s$ by a function of the form $s+h_2(u,v)$, we may
assume
\[
f_4\in V_1\coloneqq\bC u^4\oplus\bC v^4.
\]
Set
\[
\bA_{6,\mathrm{node}}\coloneqq
V_1\oplus H^0(\bP^1,\cO(6)),
\]
and let $\bG_m$ act by
\[
\lambda\cdot(f_4,f_6)
=(\lambda f_4,\lambda^2f_6).
\]
Let $\bP_{6,\mathrm{node}}$ be the coarse
moduli space of
$[(\bA_{6,\mathrm{node}}\setminus\{0\})/\bG_m]$.

The action
\[
t\cdot(f_4(u,v),f_6(u,v))
=
(f_4(tu,t^{-1}v),f_6(tu,t^{-1}v))
\]
commutes with the preceding action and descends
to $\bP_{6,\mathrm{node}}$.
We say that $[D]\in\bP_{6,\mathrm{node}}$ is
\emph{GIT-(poly/semi)stable} if it is
GIT-(poly/semi)stable for this action with
respect to $\cO_{\bP_{6,\mathrm{node}}}(1)$,
viewed as a $\bQ$-line bundle.
\end{defn}

\begin{theorem}\label{thm:k=gitp112-2}
Let $D\subseteq\bP(1,1,2)$ be represented by a
point $[D]\in\bP_{6,\mathrm{node}}$.  If $(\bP(1,1,2), cD)$ is K-(poly/semi)stable for some $c \in (c_2, \frac{2}{3})$, then $D$ is GIT-(poly/semi)stable in the sense of \Cref{defn:gitp112-2}. For $0<\epsilon\ll1$, the pair
$(\bP(1,1,2),(c_2+\epsilon)D)$ is
K-(poly/semi)stable if and only if $[D]$
is GIT-(poly/semi)stable in the sense of
\Cref{defn:gitp112-2}.  Furthermore, all such curves admit an isotrivial specialization to $(\bP(1,1,2), cC'_{\mathrm{dc}})$.
\end{theorem}

\begin{proof}
The proof follows the argument of
\cite[Theorem 5.8]{ADL24}.
\end{proof}

\begin{proof}[Proof of
\Cref{Second wall crossing}\textup{(3)}]
To describe the wall-crossing morphisms, first observe that \cite[Theorem 4.7]{ADL21} implies that any K-polystable pair $(X,cD)$ with $c \in (0, c_2 + \epsilon)$ has $X = \bP^1 \times \bP^1$ or $X = \bP(1,1,2)$.  Furthermore, by Theorem \ref{Second wall crossing}(2) and interpolation, $(\bP^1 \times \bP^1, c C_{\mathrm{dc}})$ is K-unstable for all $c > c_2$.  By Proposition \ref{curves meet the vertex point in a cusp}, all curves $D$ in K-polystable pairs $(\bP(1,1,2), cD)$ for $c \in [c_2, c_2 + \epsilon)$ are at worst nodal at the vertex and hence by Theorems \ref{thm:k=gitp112-1} and \ref{thm:k=gitp112-2}, it follows that the only new polystable point at $c = c_2$ is $(\bP(1,1,2), cC'_{\mathrm{dc}})$. 

Denote by $E_2^{\pm}$ the preimages of the point $[(\bP(1,1,2), cC'_{\mathrm{dc}})]$ under the morphisms $\phi_{c_2}^{\pm}$.   From the local VGIT description and [DH98, Theorem 0.2.5] or [Tha96, Corollary 1.13], it follows that $E_2^{\pm}$ are irreducible and we have $\dim E^-_2 + \dim E^+_2 + 1 = \dim \ove{\fM}^{\K}_4(c_2) = 9$.  The computation in Lemma \ref{de-double conic} proves that the strict transform of $\Sigma_9$ is contained in $E^-_2$ so $\dim E_2^{-} \ge 1$.  As the generic curve $D \subset \bP(1,1,2)$ that is nodal at the vertex has $\lct D > \frac{2}{3}$, the pair $(\bP(1,1,2), c D)$ is K-stable by interpolation for all $c \in (c_2, \frac{2}{3})$, so by Theorem \ref{thm:k=gitp112-2}, $[D]$ is GIT stable.  It then follows from Theorem \ref{thm:k=gitp112-2} that $\dim E_2^+ \ge 7$.  As $\dim E^-_2 + \dim E^+_2 + 1  = 9$, we conclude $E_2^-$ is the strict transform of $ \Sigma_9$ and $E_2^+ $ is the GIT quotient of $ \bP_{6,\mathrm{node}}$.  Thus, the wall crossing description is as claimed. 
\end{proof}

\subsection{Separating $A_7$-singularities}

A $(3,3)$-curve $C_{\mathrm{sep},A_7}$ with a
separating $A_7$-singularity is the union of a
$(1,1)$-curve and a $(2,2)$-curve meeting at a
single point.
By \cite{Fed12}, the corresponding locus
$\Sigma_4\subset\ove{\fM}^{\GIT}$ is
two-dimensional.

After a projective transformation,
$C_{\mathrm{sep},A_7}$ can be written in
affine coordinates as
\[
(x+z)\bigl(x+z+(x+z)(a_0x+a_1xz)+2x^2z^2\bigr)=0,
\]
where $(a_0,a_1)\in\bA^2$ and the singularity
is at $p\coloneqq(0,0)$.
Choosing the jet coordinate
$x'=x+z+z^4$ at $p$, the equation becomes
\[
(x')^2-z^8+(\text{higher order terms})=0,
\]
where $(z,x')$ has weight $(1,4)$.
Let
\[
\pi:X\coloneqq
\Bl_{(1,4)}(\bP^1\times\bP^1)
\longrightarrow\bP^1\times\bP^1
\]
be the corresponding weighted blowup, and
let $F$ be its exceptional divisor.

Set $c_4=\tfrac12$, and let
\[
C'_{\mathrm{sep},A_7}
\coloneqq V(u^2s^2+sv^4)\subset\bP(1,1,2).
\]

\begin{theorem}\label{thm:the-wall corrsponds-to-sep-A7}
Keep the notation introduced above. Then
\begin{enumerate}
\item[\textup{(1)}]
$(\bP^1\times\bP^1,cC_{\mathrm{sep},A_7})$
is K-semistable \textup{(}resp.\ K-polystable\textup{)}
if and only if $c\leq c_4$
\textup{(}resp.\ $c<c_4$\textup{)}.

\item[\textup{(2)}]
The K-polystable degeneration of
$(\bP^1\times\bP^1,c_4C_{\mathrm{sep},A_7})$
is isomorphic to
$(\bP(1,1,2),c_4C'_{\mathrm{sep},A_7})$.

\item[\textup{(3)}]
In a neighborhood of the point $[(\bP(1,1,2),c_4C'_{\mathrm{sep},A_7})] \in \ove{\fM}^{\K}_4(c_4) $ the wall crossing is given by a flip of the
strict transform of $\Sigma_4$.
The associated exceptional loci
$E^+_{4,\mathrm{sep}}$ and
$E^-_{4,\mathrm{sep}}$ have dimensions $6$
and $2$, respectively.
The former parametrizes K-polystable pairs
$(\bP(1,1,2),(c_4+\epsilon)D)$ with
\[
D=V\bigl(u^2s^2+sv^4+g_6(u,v)\bigr),
\]
where $g_6\neq0$ is a binary sextic.
\end{enumerate}
\end{theorem}

\begin{proof}
We first show that
$(\bP^1\times\bP^1,cC_{\mathrm{sep},A_7})$
is K-unstable for $c>c_4$. We have
\[
A_{(\bP^1\times\bP^1,cC_{\mathrm{sep},A_7})}(F)
=5-8c.
\]
A direct computation gives
\[
\vol\bigl(-\pi^*K_{\bP^1 \times \bP^1}-tF\bigr)
=
\begin{cases}
8-\tfrac14t^2,
&0\leq t\leq4,\\[2pt]
\tfrac14(8-t)^2,
&4\leq t\leq8,
\end{cases}
\]
and hence
\[
S_{(\bP^1\times\bP^1,cC_{\mathrm{sep},A_7})}(F)
=4-6c.
\]
Therefore,
\[
\beta_{(\bP^1\times\bP^1,cC_{\mathrm{sep},A_7})}(F)
=1-2c<0
\qquad\text{for }c>c_4.
\]

\begin{lem}\label{de-sep A_7}
The pair
$(\bP^1\times\bP^1,cC_{\mathrm{sep},A_7})$
admits a special degeneration to
$(\bP(1,1,2),cC'_{\mathrm{sep},A_7})$.
\end{lem}

\begin{proof}
Let $\sX_0=(\bP^1\times\bP^1)\times\bA^1_t$,
and let $\pi:\sX\to\sX_0$ be the
$(1,1,4)$-weighted blowup in the local
coordinates $(t,z,x')$.
Its central fiber is $S\cup X$, where
$S\simeq\bP(1,1,4)$ is the exceptional divisor
and $X$ is the strict transform of the
original central fiber.

Let $Q=V(x+z)$, and let $\overline Q\subset X$
be its strict transform. Since
$x+z=x'-z^4$, we have
\[
\pi_0^*Q=\overline Q+4F.
\]
The curve $\overline Q$ meets $F$ at a smooth
point of both $S$ and $X$, and
\[
\overline Q^2=-2,\qquad
F\cdot\overline Q=1.
\]
Moreover,
\[
N_{\overline Q/\sX}
\simeq\cO_{\bP^1}(-2)\oplus\cO_{\bP^1}(-1).
\]

We perform the birational transformations
\[
\begin{tikzcd}
&\sX
  \arrow[ld,"\pi"']
  \arrow[rd,"g"]
  \arrow[rr,dashed,"f"]
&&\sX^+
  \arrow[ld,"h"']
  \arrow[rd,"\psi"]
\\
\sX_0&&\sY&&\sZ,
\end{tikzcd}
\]
where $g$ contracts $\overline Q$, $f$ is the
corresponding flip for the relative
$(-K_{\sX}-S)$-MMP, and $\psi$ contracts
the strict transform $X'$ of $X$ to a point.

The induced morphism $\hat S\to S$ is the
$(2,1)$-weighted blowup at
$\overline Q\cap F$.
The flipped curve meets the strict transform
$\hat F$ of $F$ at an $A_1$-singularity.
The restriction of $\psi$ to $\hat S$
contracts $\hat F$, and the resulting surface
is $S'\simeq\bP(1,1,2)$.
Tracking the strict transform of
$C_{\mathrm{sep},A_7}\times\bA^1$
gives the limiting curve
\[
C'_{\mathrm{sep},A_7}=V(u^2s^2+sv^4).
\]
\end{proof}

The final central fiber is the birational
transform of $S$, and
\[
\ord_S(t)=1,\qquad
\ord_S|_{K(\bP^1\times\bP^1)}=\ord_F.
\]
Thus this special test configuration is
induced by $\ord_F$.

The pair
$(\bP(1,1,2),c_4C'_{\mathrm{sep},A_7})$
is a $\mathbb T$-pair of complexity one
under the action
\[
[u:v:s]\longmapsto[u:\lambda v:\lambda^4s].
\]
There are no horizontal divisors, and the
vertical prime divisors are
\[
V(u),\qquad V(v),\qquad V(s),\qquad
V(u^2s-av^4),\qquad a\in\bC^*.
\]
At $c=c_4$, their $\beta$-invariants are
positive, and
\[
\operatorname{Fut}_{
(\bP(1,1,2),c_4C'_{\mathrm{sep},A_7})}(\lambda)
=1-2c_4=0.
\]
Hence the pair is K-polystable by
\Cref{thm:complexity-one}.
Openness of K-semistability implies that
$(\bP^1\times\bP^1,c_4C_{\mathrm{sep},A_7})$
is K-semistable.
Together with
\Cref{thm:first-wall}\textup{(1)} and
interpolation of K-stability, this gives
K-polystability for $0<c<c_4$.

To describe the wall-crossing morphisms, first observe that any K-polystable pair $(X,(c_4 \pm \epsilon) D)$ in a neighborhood of $(\bP(1,1,2),c_4C'_{\mathrm{sep},A_7})$ has $X = \bP^1 \times \bP^1$ or $X = \bP(1,1,2)$ as these are the only $\mathbb{Q}$-Gorenstein deformations of the surface $\bP(1,1,2)$.  Furthermore, by Theorem \ref{thm:the-wall corrsponds-to-sep-A7}(2) and interpolation, $(\bP^1 \times \bP^1, c C_{\mathrm{sep},A_7})$ is K-unstable for all $c > c_4$.  By Proposition \ref{curves meet the vertex point in a cusp}, all curves $D$ in K-polystable pairs $(\bP(1,1,2), cD)$ for $c \in [c_4, c_4 + \epsilon)$ are at worst cuspidal at the vertex and hence by Theorems \ref{thm:k=gitp112-1}, \ref{thm:k=gitp112-2}, and Proposition \ref{singularities of curves meet the vertex point}, it follows that the only new polystable point $(X,cD)$ with $X = \bP(1,1,2)$ at $c = c_4$ is $(\bP(1,1,2), c_4 C'_{\mathrm{sep},A_7})$. 

Denote by $E_{4,\mathrm{sep}}^{\pm}$ the preimages of the point $[(\bP(1,1,2), c_4 C'_{\mathrm{sep},A_7})]$ under the morphisms $\phi_{c_4}^{\pm}$. From the local VGIT description and [DH98, Theorem 0.2.5] or [Tha96, Corollary 1.13], it follows that $E_{4,\mathrm{sep}}^{\pm}$ are irreducible and we have $\dim E_{4,\mathrm{sep}}^{-} + \dim E_{4,\mathrm{sep}}^{+} + 1 = \dim \ove{\fM}^{\K}_4(c_4) = 9$.  The computation in Lemma \ref{de-sep A_7} proves that the strict transform of $\Sigma_4$ is contained in $E_{4,\mathrm{sep}}^{-}$ so $\dim E_{4,\mathrm{sep}}^{-} \ge 2$.  As the generic curve $D \subset \bP(1,1,2)$ that is cuspidal at the vertex has $\lct D > \frac{2}{3}$, the pair $(\bP(1,1,2), c D)$ is K-stable by interpolation for all $c \in (c_4, \frac{2}{3})$.  Up to coordinate change, pairs with $D$ cuspidal at the vertex are given by an equation of the form
\[
u^2s^2+sv^4+g_6(u,v)=0,
\qquad g_6\neq0.
\]
The binary sextic contributes $7$ parameters,
and the remaining scaling reduces the
dimension of this locus to $6$.

It then follows that $\dim E_{4,\mathrm{sep}}^{+} \ge 6$.  As $\dim E_{4,\mathrm{sep}}^{-} + \dim E_{4,\mathrm{sep}}^{+} + 1  = 9$, we conclude $E_{4,\mathrm{sep}}^{-} = \Sigma_4$ and $E_{4,\mathrm{sep}}^{+}$ parametrizes the K-polystable pairs $(\bP(1,1,2), (c_4 + \epsilon)D)$ such that $D$ is cuspidal at the vertex.  Thus, the wall crossing description is as claimed. 
\end{proof}

\begin{prop}\label{singularities of curves meet the vertex point}
Let $(\bP(1,1,2),(c_4+\epsilon)D)$ be
K-semistable, where $0<\epsilon\ll1$.
If $D$ passes through the vertex and does
not have a node there, then it has a cusp
there. A general such curve is smooth away
from the vertex.

If $D$ is cuspidal at the vertex, then $(\bP(1,1,2), c_4D)$ isotrivially specializes to $(\bP(1,1,2), c_4 C'_{\mathrm{sep},A_7})$.
\end{prop}

\begin{proof}
The first statement follows from
\cite[Proposition 4.12]{CMJL14}
and the K-moduli and VGIT quotient comparison in \cite{ADL21} and the results above.  

The second statement follows from the fact that, up to coordinate change, a curve cuspidal at the vertex is given by an equation of the form
\[
u^2s^2+sv^4+g_6(u,v)=0,
\qquad g_6\neq0.
\]
\end{proof}

\section{The wall crossing involving $X_{10}$}
\label{sec:walls of X10}

By \Cref{strata-for-GIT}, the locus $\Sigma_2\subseteq \ove{\fM}^{\GIT}$
parametrizing the orbits of $(3,3)$-curves with a separating $A_9$-singularity consists of a single point. Set $c_3=\tfrac{16}{33}$, and let $[x:y:z:w]$ be the weighted homogeneous coordinates of $\bP(1,1,5,9)$, and
\[
X_{10} \ \coloneqq\ 
V(xw-y^{10}-z^2)\ \subseteq \ \bP(1,1,5,9),
\qquad
C'_{A_9}\ \coloneqq\  V(w)\ \subseteq\ X_{10}.
\] The surface $X_{10}$ is a partial smoothing of $\bP(1,2,9)$ obtained by smoothing its $A_1$-singularity. Indeed, the family
\[
    V(xw-z^2-t y^{10})\ \subseteq\ \bP(1,1,5,9)\times\bA^1_t
\]
has central fiber isomorphic to $\bP(1,2,9)$ and all other fibers isomorphic to $X_{10}$.

\begin{theorem}\label{Main Theorem of X10}
Let $C_{A_9}$ be the $(3,3)$-curve with a separating
$A_9$-singularity. Then
\begin{enumerate}
\item[\textup{(1)}]
The wall-crossing morphism
\[
\phi^-_{c_3}\ \colon\ 
\ove{\fM}^{\K}_4(c_3-\epsilon)
\ \longrightarrow\ 
\ove{\fM}^{\K}_4(c_3)
\]
locally over $[(X_{10},c_3C'_{A_9})]$ is an isomorphism, sending $[(\bP^1\times\bP^1,(c_3-\epsilon)C_{A_9})]$
to $[(X_{10},c_3C'_{A_9})]$. 
\item[\textup{(2)}]
the wall-crossing morphism
\[
\phi^+_{c_3} \ \colon\ 
\ove{\fM}^{\K}_4(c_3+\epsilon)
\ \longrightarrow\ 
\ove{\fM}^{\K}_4(c_3)
\]
locally over $[(X_{10},c_3C'_{A_9})]$ is a divisorial contraction.
Its exceptional divisor parametrizes pairs $(X_{10},(c_3+\epsilon)C')$ with
\[
C'\ =\ V(w-g_9(x,y))\cap X_{10},
\]
where $g_9\neq0$ is homogeneous of degree $9$;
\item[\textup{(3)}]
for any curve $C'$ of weighted degree $9$ passing through the unique singular point of $X_{10}$, $(X_{10},cC')$ is K-unstable for any $0<c<\tfrac23$.
\end{enumerate}
\end{theorem}

We prove the theorem in the following two subsections.

\subsection{K-polystable replacement of the separating
$A_9$-curve}

After a projective transformation, $C_{A_9}$ can be
written in affine coordinates as
\[
r\bigl(z^5+r(z^3+2z^2+xz+4x+4z+8)\bigr)=0,
\qquad
r\coloneqq x+2z-z^2,
\]
with its singularity at $p\coloneqq(0,0)$.
Choosing the jet coordinate
$x'=x+2z-z^2+\tfrac{z^5}{16}$ at $p$, the equation
becomes
\[
8(x')^2-\tfrac{1}{32}z^{10}
+(\text{higher order terms})=0,
\]
where $(z,x')$ has weight $(1,5)$.
Let
\[
\pi\ :\ X\coloneqq
\Bl_{(1,5)}(\bP^1\times\bP^1)
\ \longrightarrow\ \bP^1\times\bP^1
\]
be the corresponding weighted blowup at $p$, and
let $F$ be its exceptional divisor.

\begin{prop}\label{The wall corresponds to sep A9}
Keep the notation introduced above. Then
\begin{enumerate}
\item[\textup{(1)}]
$(\bP^1\times\bP^1,cC_{A_9})$ is K-semistable
\textup{(}resp.\ K-polystable\textup{)} if and only if
$c\leq c_3$ \textup{(}resp.\ $c<c_3$\textup{)};
\item[\textup{(2)}]
the K-polystable degeneration of
$(\bP^1\times\bP^1,c_3C_{A_9})$ is isomorphic to
$(X_{10},c_3C'_{A_9})$.
\end{enumerate}
\end{prop}

\begin{proof}
We first show that $(\bP^1\times\bP^1,cC_{A_9})$
is K-unstable for $c>c_3$. We have $A_{(\bP^1\times\bP^1,cC_{A_9})}(F)=6-10c$. A direct computation gives
\[
\vol\bigl(-\pi^*K_{\bP^1 \times \bP^1}-tF\bigr)
=
\begin{cases}
8-\tfrac15t^2,
&0\leq t\leq6,\\[2pt]
\tfrac95(\tfrac{20}{3}-t)^2,
&6\leq t\leq\tfrac{20}{3},
\end{cases}
\]
and hence for any $c>c_3$
\[
S_{(\bP^1\times\bP^1,cC_{A_9})}(F)
\ =\ \tfrac{38}{9}-\tfrac{19}{3}c,\qquad 
\beta_{(\bP^1\times\bP^1,cC_{A_9})}(F)
\ =\ \tfrac{16}{9}-\tfrac{11}{3}c\ <\ 0
\qquad\text{for }c>c_3.
\]

\begin{lem}\label{de-Sep A9}
The pair $(\bP^1\times\bP^1,cC_{A_9})$ admits a special
degeneration to $(X_{10},cC'_{A_9})$.
\end{lem}

\begin{proof}
Let $\phi:\sX\to(\bP^1\times\bP^1)\times\bA^1_t$ be the $(1,1,5)$-weighted
blowup in the local coordinates $(t,z,x')$.
Its central fiber is $S\cup X$, where
$S\simeq\bP(1,1,5)$ is the exceptional divisor and
$X$ is the strict transform of the original central
fiber.

Write $C_{A_9}=Q_0\cup Q_1$, where $Q_0$ and $Q_1$
have bidegrees $(1,2)$ and $(2,1)$, respectively,
and let $\tQ_i\subseteq X$ be their strict
transforms. The curves $\tQ_0$ and $\tQ_1$ are
disjoint $(-1)$-curves on $X$. Each meets
$F=S\cap X$ transversely at a smooth point $p_i$,
and
\[
N_{\tQ_i/\sX}
\ \simeq\ \cO_{\bP^1}(-1)\oplus\cO_{\bP^1}(-1).
\] We run a $(-K_{\sX}-S)$-MMP over $\bA^1$:
\[
\begin{tikzcd}
&\sX
  \arrow[ld,"\phi"']
  \arrow[rd,"g"]
  \arrow[rr,dashed,"f"]
&&\sX^+
  \arrow[ld,"h"']
  \arrow[rd,"\psi"]
\\
(\bP^1\times\bP^1)\times\bA^1&&\sY&&\sZ.
\end{tikzcd}
\]
The first contraction $g$ contracts the curves $\tQ_0$ and $\tQ_1$: for $i=0,1$, we have
\[
    (-K_{\sX}-S)\cdot\tQ_i\ =\ -1,
    \qquad
    K_{\sX}\cdot\tQ_i=0.
\]
Thus the corresponding $(-K_{\sX}-S)$-flip
$f\colon\sX\dashrightarrow\sX^+$
is the simultaneous Atiyah flop of these two curves, and its restriction to $S$ gives
\[
h:\hat S=\Bl_{p_0,p_1}S\longrightarrow S.
\] The next step $\psi$ contracts the strict transform $X'$ of $X$ to a point, and its restriction to $\hat S$ contracts
the strict transform $\hat F$ of $F$. On the central fibers, these maps give
\[
\begin{tikzcd}
&S\cup X
  \arrow[ld,"\pi"']
  \arrow[rd,"g"]
  \arrow[rr,dashed,"f"]
&&\hat S\cup X'
  \arrow[ld,"h"']
  \arrow[rd,"\psi"]
\\
\bP^1\times\bP^1&&S\cup X'&&S'.
\end{tikzcd}
\] After a change of coordinates, we can write
$S=\bP(1,1,5)_{[x_1:x_2:x_3]}$ so that
\[
F=V(x_1),\qquad
\{p_0,p_1\}\ =\ V(x_1,x_3^2+x_2^{10}).
\]
Let $E_0,E_1$ be the $h$-exceptional curves, and let $D$ be $E_0+E_1+\frac{10}{9}\hat F$. Since one has the intersection numbers \[
    D^2\ =\ \tfrac{2}{9},\qquad
    (D\cdot\hat F)\ =\ 0,\qquad
    (D\cdot E_i)\ =\ \tfrac{1}{9},
\] $D$ is nef and big. Moreover, since $(\hat S,\frac13\hat F)$ is klt and $K_{\hat S}+\frac13\hat F\sim_{\bQ}-6D$, $D$ is
semiample by the base-point-free theorem. The $D$-ample model morphism contracts precisely $\hat F$, and hence it is exactly $\psi|_{\hat S}:\hat S\to S'$. Writing $H=h^*F$, we have
\[
    \lfloor mD\rfloor
    \ =\ mH+\left\lfloor\tfrac{m}{9}\right\rfloor\hat F.
\]
Thus a section of $\cO_{\hat S}(\lfloor mD\rfloor)$
can be represented by $P/x_1^k$, where
$0\leq k\leq\lfloor m/9\rfloor$, $P$ is weighted homogeneous of degree $m+k$, and $P$ vanishes to order at least $k$ at both $p_0$ and $p_1$. As $\{p_0,p_1\}=V(x_1,x_3^2+x_2^{10})$, $P$ can be written as
\[
    P\ =\ \sum_{j=0}^k x_1^{k-j}(x_3^2+x_2^{10})^jP_j,
    \qquad \deg P_j=m-9j.
\]
Consequently, the section ring of $D$ is generated by
\[
    x=x_1,\qquad y=x_2,\qquad z=x_3,\qquad
    w=\frac{x_3^2+x_2^{10}}{x_1},
\]
of degrees $1,1,5,9$, respectively, with the single relation $xw=z^2+y^{10}$, and it follows that
\[
    S'\ \simeq\ 
    \Proj
    \frac{\bC[x,y,z,w]}{(xw-y^{10}-z^2)}
    \ =\ X_{10}.
\] The strict transform of $C_{A_9}\times\bA^1$
meets $S$ along $V(x_3^2+x_2^{10})$.
Tracking this curve through the flop and the final contraction gives $V(w)=C'_{A_9}$ on $X_{10}$.
\end{proof}

The resulting special test configuration $\sZ\to \bA^1$ is induced by $\ord_F$. The pair $(X_{10},c_3C'_{A_9})$ is a
$\mathbb T$-pair of complexity one under the action
\[
[x:y:z:w]\ \longmapsto\ 
[\lambda^{10}x:\lambda y:\lambda^5z:w].
\]
On $X_{10}$, there are no horizontal divisors, and the vertical prime divisors are
\[
V(y),\qquad
V(z-ay^5),\qquad V(x,z-iy^5),\qquad V(x,z+iy^5),\qquad
V(w,z-iy^5),\qquad V(w,z+iy^5),
\]
where $a\in\bC\setminus\{i,-i\}$. At $c=c_3$, their $\beta$-invariants are all positive, and
\[
\Fut_{(X_{10},c_3C'_{A_9})}(\lambda)\ =\ \tfrac{16-33c_3}{9}\ =\ 0.
\]
Thus $(X_{10},c_3C'_{A_9})$ is K-polystable by
\Cref{thm:complexity-one}. The openness of K-semistability implies that
$(\bP^1\times\bP^1,c_3C_{A_9})$ is K-semistable, and
the interpolation of K-stability then shows that
$(\bP^1\times\bP^1,cC_{A_9})$ is K-stable for $0<c<c_3$.
\end{proof}

\subsection{K-stability of curves on $X_{10}$}

We next study curves of weighted degree $9$ on $X_{10}$.

\begin{lem}\label{lem:curves-on-X10}
Let $C'\subseteq X_{10}$ be a curve of weighted degree $9$. 
\begin{enumerate}
    \item[\textup{(1)}] If $(X_{10},cC')$ is K-semistable, then $c\geq c_3$.
    \item[\textup{(2)}] If $C'$ passes through the unique singular point of $X_{10}$, then $(X_{10},cC')$ is K-unstable for every $0<c<\tfrac23$.
\end{enumerate} 
\end{lem}

\begin{proof}
Let $\mu\colon\hat S\to X_{10}$ be the contraction
constructed above, and set $E=\hat F$. Its image is
the unique singular point $p=[0:0:0:1]$ of $X_{10}$.
We have
\[
    E^2=-\tfrac{9}{5},
    \qquad
    K_{\hat S}=\mu^*K_{X_{10}}-\tfrac{1}{3}E, \qquad A_{X_{10}}(E)=\tfrac{2}{3}.
\] Set $T=V_{X_{10}}(x)$. Its strict transform is
$\overline T=E_0+E_1$, so
\[
    \mu^*T\ =\ E_0+E_1+\tfrac{10}{9}E,
    \qquad
    \ove{T}^2=-2.
\]
Using $-\mu^*K_{X_{10}}
    \sim_{\bQ}6(E_0+E_1)+\tfrac{20}{3}E$, the Zariski decomposition gives
\[
\vol\bigl(-\mu^*K_{X_{10}}-tE\bigr)
\ =\ 
\begin{cases}
8-\frac95t^2,
&0\leq t\leq\frac23,\\[2pt]
\frac15\left(\frac{20}{3}-t\right)^2,
&\frac23\leq t\leq\frac{20}{3},
\end{cases}\qquad S_{(X_{10},cC')}(E)
    \ =\ \tfrac{22}{9}-\tfrac{11}{3}c.
\]
Consequently,
\[
    \beta_{(X_{10},cC')}(E)
    \ =\ \tfrac{2}{3}-c\,\ord_E(C')
      -\left(\tfrac{22}{9}-\tfrac{11}{3}c\right)\ \leq\ -\tfrac{16}{9}+\tfrac{11}{3}c.
\] This is negative for $c<c_3=\frac{16}{33}$,
proving the first assertion. If, moreover, $C'$ passes through $p$, then $\ord_E(C')\geq1$, and hence
\[
    \beta_{(X_{10},cC')}(E)
    \ \leq\ -\tfrac{16}{9}+\tfrac{8}{3}c\ <\ 0
\]
for any $0<c<\frac23$, which proves the second assertion.
\end{proof}

\begin{proof}[Proof of \Cref{Main Theorem of X10}]
By \Cref{The wall corresponds to sep A9}, the pair
$(\bP^1\times\bP^1,c_3C_{A_9})$ has K-polystable replacement $(X_{10},c_3C'_{A_9})$.
By \Cref{lem:curves-on-X10}, every curve $C'$ on $X_{10}$ of weighted degree $9$ such that $\big(X_{10},(c_3+\epsilon)C'\big)$ is K-semistable for $0<\epsilon\ll1$ avoids the singular point of $X_{10}$. Up to an automorphism of $X_{10}$, such a curve can be written as
\[
    C_{\mathbf b}\ =\ 
    V\Big(w-\sum_{i=1}^{9}b_i x^iy^{9-i}\Big)\cap X_{10}.
\]
The residual $\bG_m$-action on the parameters, induced by
\[
[x:y:z:w]\mapsto[\lambda x:y:z:\lambda^{-1}w],
\]
has weights $2,3,\ldots,10$. Since every $C_{\mathbf b}$ admits an equivariant degeneration to $C'_{A_9}$, the pair $(X_{10},c_3C_{\mathbf b})$ is K-semistable. Thus the nonzero parameters give an $8$-dimensional locus contracted by $\phi^+_{c_3}$ to $[(X_{10},c_3C'_{A_9})]$. Let $E_3\subseteq \ove{\fM}^\K_4(c_3)$ be the connected component of the wall-crossing exceptional locus containing $[(X_{10},c_3C'_{A_9})]$. Since $\Aut^0(X_{10},c_3C'_{A_9})=\bG_m$, it follows from \cite[Theorem 0.2.5]{DH98} that $\dim E_3=0$. In particular, $E_3=\{[(X_{10},c_3C'_{A_9})]\}$ and $(\phi^-_{c_3})^{-1}(E_3) = \{[(\bP^1 \times \bP^1,c_3C_{A_9})]\}$. By Zariski's Main Theorem, $\phi^-_{c_3}$ locally over $[(X_{10},c_3C'_{A_9})]$ is an isomorphism, and hence $E_3^+$ is the closure of the locus parametrizing $(X_{10},(c_3+\epsilon)C_{\mathbf b})$. This proves \textup{(2)}. Together with \Cref{strata-for-GIT,lem:curves-on-X10}, this identifies the unique replacement at $c_3$ locally near $[(X_{10},c_3C'_{A_9})]$.

The final assertion is \Cref{lem:curves-on-X10}(2).
\end{proof}

\begin{remark}
In the above normal form, the $y^9$ coefficient has been eliminated by an automorphism of $X_{10}$. The resulting parameter quotient is $\bP(2,3,\ldots,10)$, which describes $E_3^+$. However, its strict transform is modified by subsequent wall-crossing, unlike the unigonal locus in \cite[Section 5]{ADL22}.
\end{remark}

\section{Wall-crossings involving $\mathbb P(1,2,9)$} \label{sec: walls of P(1,2,9)}

In this section, we treat $(3,3)$-curves with a non-separating $A_{k}$-singularity for $6 \le k \le 8$. We will see that the wall-crossing morphisms corresponding to these curves are all flips, and the new pairs appearing across the wall are curves on $\bP(1,2,9)$.

\subsection{$A_{8}$-singularity}
By \Cref{strata-for-GIT}, the locus $\Sigma_3 \subseteq \ove{\fM}^{\GIT}$ parametrizing $(3,3)$-curves with an $A_{8}$-singularity is one-dimensional, and its closure contains the stratum $\Sigma_2$ parametrizing the unique orbit of the $(3,3)$-curves with a separating $A_{9}$-singularity. After a projective transformation, in the affine coordinates a $(3,3)$-curve $C_{A_8}$ with an $A_8$-singularity is defined by the polynomial \[
   (x + 2z - z^{2})\big(z^{5} + (x + 2z - z^{2})(z^{3} + 2z^{2} + xz + 4x + 4z +8)\big) + a_{0}\big(xz + 2(x + 2z)\big)^3 ,\]
 where $a_{0} \in \bA^1 \setminus \{0\}$. Choosing the jet coordinate $x' = x + 2z - z^{2} + \tfrac{z^{5}}{16}$ at $p\coloneqq (0,0)$, the equation of $C_{A_8}$ in $(x',z)$ becomes
\begin{equation*}
    8(x')^{2} +a_{0}z^{9} + (\text{higher order terms})=0,
\end{equation*}
where $(z,x')$ has weight $(2,9)$. Let \[\pi:X\coloneqq \Bl_{(2,9)}\bP^1\times\bP^1\ \longrightarrow \  \bP^1\times\bP^1\] be the corresponding $(2,9)$-weighted blowup at $p$, and $F$ be the exceptional divisor.  Let $H \in |\cO_{\bP^1 \times \bP^1}(1,1)|$ denote the class of the hyperplane section.

\begin{lem}\label{lem:bpfness-1}
    The divisor $3\pi^*H-18F$ is base-point free, and induces a fibration over $F$.
\end{lem}

\begin{proof}
Set
\[
r\coloneqq x+2z-z^2,\qquad q\coloneqq xz+2x+4z,\qquad
s\coloneqq x^2z+4x^2+12xz+8x+16z.
\]
The pencil generated by $rs$ and $q^3$ has its unique
base point at $p$: $V(q)$ meets each of $V(r)$ and
$V(s)$ only at $p$, with intersection multiplicity $3$. In the local coordinates $(z,x')$ of weights $(2,9)$, we have
\[
rs=8(x')^2+ (\text{higher order terms}),\qquad
q^3=z^9+ (\text{higher order terms}).
\]
Thus both generators have weighted order $18$, and their strict
transforms have no common zero on $F\simeq\bP(2,9)$.
They therefore generate
$3\pi^*H-18F$ and define a morphism
$g:X\to\bP^1$.
Its restriction to $F$ is the isomorphism
\[
[x':z]\ \longmapsto\ [8(x')^2:a_0z^9],
\]
hence $g$ has connected fibers. Identifying $\bP^1$ with $F$
gives the desired fibration.
\end{proof}

Set $c_4=\frac12$, and let $C'_{A_8}\subset\bP(1,2,9)$ be the curve defined by $v^9+s^2=0$. Up to automorphisms of $\bP(1,2,9)$, $C'_{A_8}$ is the unique curve of weighted degree $18$ whose only singularity is of type $A_8$ and lies in the smooth locus of $\bP(1,2,9)$.

\begin{theorem} \label{the wall corresponds to A8} Keep the notation introduced above. Then 
\begin{enumerate}
    \item[\textup{(1)}] $(\bP^{1} \times \bP^{1}, cC_{A_8})$ is K-semistable \textup{(}resp.\ K-polystable\textup{)} if and only if $c \leq c_4$ \textup{(}resp.\ $c < c_4$\textup{)};
    \item[\textup{(2)}] the K-polystable degeneration of $(\bP^{1} \times \bP^{1}, c_4C_{A_8})$ is isomorphic to $\big(\bP(1,2,9), c_4C'_{A_8}\big)$;
    \item[\textup{(3)}] 
$\phi_{c_4}\colon \ove{\fM}^\K_4(c_4-\epsilon)\dashrightarrow\ove{\fM}^\K_4(c_4+\epsilon)$ locally over $[\big(\bP(1,2,9), c_4C'_{A_8}\big)]$ is a flip, with the exceptional loci of $\phi_{c_4}^{+}$ and $\phi_{c_4}^{-}$ having dimensions $7$ and $1$, respectively.
\end{enumerate} 
\end{theorem}
\begin{proof}
    We first show that $(\bP^{1} \times \bP^{1}, cC_{A_8})$ is K-unstable for any $c > \tfrac{1}{2}$. It is straightforward that $A_{(\bP^{1} \times \bP^{1}, cC_{A_8})}(F) = 11 - 18c$. By \Cref{lem:bpfness-1}, $2\pi^{*}H - tF$ is nef if and only if $0 \leq t \leq 12$, and a direct computation then gives  
$S_{(\bP^{1} \times \bP^{1}, cC_{A_8})}(F) = 8 - 12c$. Then for any $c > \tfrac{1}{2}$, one has
\[\beta_{(\bP^{1} \times \bP^{1}, cC_{A_8})}(F) \ = \ A_{(\bP^{1} \times \bP^{1}, cC_{A_8})}(F) \ - \ S_{(\bP^{1} \times \bP^{1}, cC_{A_8})}(F) \ = \ 3 \ - \ 6c\ <\ 0.\]

\begin{lem}\label{de-A_8}
    The pair $(\bP^1 \times \bP^1,cC_{A_8})$ admits a special degeneration to $(\bP(1,2,9),cC'_{A_8})$.
\end{lem}
\begin{proof}
    Let $\sX \coloneqq (\bP^1\times\bP^1)\times\bA^1_t$, and let $\phi:\sX'\rightarrow \sX$ be the $(1,2,9)$-weighted blowup of $\sX$ in the local coordinates $(t,z,x')$. Then the exceptional divisor $S$ is isomorphic to $\bP(1,2,9)$. We know that $\sX^{'}_0 = S \cup X$, where $X$ is the strict transform of $\sX_0$. We claim that $-K_{\sX^{'}}-S$ is $\pi$-semiample. For $t \neq 0$, we have $\sX^{'}_t \simeq \bP^1 \times \bP^1$, and $(-K_{\sX^{'}}-S)|_{\sX^{'}_t} \sim 2H_t$ is ample. For $t = 0$, $(-K_{\sX^{'}}-S)|_{S} \sim -K_S $ is ample; and $(-K_{\sX^{'}}-S)|_{X}\sim \pi^*\cO(2,2)-12F$ is semiample by \Cref{lem:bpfness-1}. Thus, the ample model of $-K_{\sX^{'}}-S$ morphism contracts the $X$ to the curve $F=X \cap S$ and is an isomorphism elsewhere. Tracking the curve $C_{A_8}\times \bA^1$ through this procedure, the limit of $C_t$ is exactly $C'_{A_8}\subseteq \bP(1,2,9)$.
\end{proof}

By construction, this special test configuration is induced by the divisorial valuation $\ord_{F}$. The pair $(\bP(1,2,9), \tfrac{1}{2}C'_{A_8})$ is a $\mathbb{T}$-pair of complexity one with the $\mathbb{G}_m$-action $\lambda:[u:v:s]\mapsto [u:\lambda^2v:\lambda^9s]$. There is a unique horizontal divisor $E\coloneqq V(u)$ on $\bP(1,2,9)$, and the vertical divisors on $\bP(1,2,9)$ are $V(v)$,  $V(s)$, and $V(s^2-av^9)$ for $a \in \bC^{*}$. It is easy to check that for $c_4= \frac{1}{2}$, the $\beta$-invariants of the vertical divisors are all positive and 
\begin{align*}
    \beta_{(\bP(1,2,9), c_4C'_{A_8})}(E) \ =\  6c_4 -3 \ =\ 0, \\
    \Fut_{(\bP(1,2,9), c_4C'_{A_8})}(\lambda) \ =\  3- 6c_4  \ =\ 0.
\end{align*}
Thus $(\bP(1,2,9), \tfrac{1}{2}C_{A_8}')$ is K-polystable by \Cref{thm:complexity-one}, and $(\bP^{1} \times \bP^{1}, c_4C_{A_8})$ is K-semistable due to the openness of K-semistability. By interpolation of K-stability, $(\bP^{1} \times \bP^{1}, cC_{A_8})$ is K-stable for any $0<c<c_4$.

\smallskip

For any curve $C'=V(s^2+v^9+u^4g_{14}(u,v))\subset\bP(1,2,9)$, there is a $\bG_m$-equivariant degeneration of $C'$ to $C'_{A_8}$. Consequently, $(\bP(1,2,9),c_4C')$ is K-semistable. If $C'$ is smooth, then $(\bP(1,2,9),\tfrac23C')$ is klt, so interpolation of K-stability implies that $(\bP(1,2,9),cC')$ is K-stable for every $c_4<c<\tfrac23$. On the other hand, the $\beta$-invariant with respect to $E$ shows that $(\bP(1,2,9),cC')$ is K-unstable for $c<c_4$. Thus the $\phi^+_{c_4}$-exceptional locus $E^+_4$ contains the closure of the locus of these smooth pairs, which has dimension $7$ by a parameter count. The proper transform of the $A_8$-locus is contained in $E^-_{4}$ and has dimension $1$. Applying \cite[Theorem 0.2.5]{DH98} to the local VGIT presentation at $[(\bP(1,2,9),c_4C'_{A_8})]$, we deduce that $E_4$ is zero-dimensional near this point. Hence, locally,
\[
    E_4\ =\ \{[(\bP(1,2,9),c_4C'_{A_8})]\},
\]
$E^-_4$ is the proper transform of $\ove{\Sigma}_3$, and $E^+_4$ is the closure of the locus described above. In particular, $\phi_{c_4}$ is a flip locally over $[(\bP(1,2,9),c_4C'_{A_8})]$.
\end{proof}

\subsection{Non-separating $A_{7}$-singularity}

After a projective transformation, a $(3,3)$-curve $C_{A_7}$
with a non-separating $A_7$-singularity at $p\coloneqq(0,0)$
can be written in affine coordinates as
\[
r\bigl(z^5+r(z^3+2z^2+xz+4x+4z+8)\bigr)
+a_0q^3+a_1q^2(x+2z)\ =\ 0,
\]
where
\[
r\coloneqq x+2z-z^2,\qquad q\coloneqq xz+2x+4z,
\]
and $(a_0,a_1)\in\bA^1\times(\bA^1\setminus\{0\})$.
Choosing the jet coordinate
$x'=x+2z-z^2+\tfrac{z^5}{16}$ at $p$, the equation becomes
\[
8(x')^2+a_1z^8+(\text{higher order terms})=0,
\]
where $(z,x')$ has weight $(1,4)$.
Let
\[
\pi:X\coloneqq\Bl_{(1,4)}(\bP^1\times\bP^1)
\ \longrightarrow\ \bP^1\times\bP^1
\]
be the corresponding weighted blowup at $p$, and let $F$
be its exceptional divisor.  Let $H \in |\cO_{\bP^1 \times \bP^1}(1,1)|$ denote the class of the hyperplane section.

\smallskip

Set $c_5=\tfrac{11}{21}$, and let
$C'_{A_7}\subseteq\bP(1,2,9)$ be the curve defined by $s^2+u^2v^8=0$.

\begin{theorem}\label{The wall corresponds to non sep A7}
Keep the notation introduced above. Then
\begin{enumerate}
\item[\textup{(1)}]
$(\bP^1\times\bP^1,cC_{A_7})$ is K-semistable
\textup{(}resp.\ K-polystable\textup{)} if and only if
$c\leq c_5$ \textup{(}resp.\ $c<c_5$\textup{)};
\item[\textup{(2)}]
the K-polystable degeneration of
$(\bP^1\times\bP^1,c_5C_{A_7})$ is isomorphic to
$(\bP(1,2,9),c_5C'_{A_7})$;
\item[\textup{(3)}]
$\phi_{c_5}\colon
\ove{\fM}^{\K}_4(c_5-\epsilon)\dashrightarrow
\ove{\fM}^{\K}_4(c_5+\epsilon)$, locally over $[(\bP(1,2,9),c_5C'_{A_7})]$, is a flip, with the
exceptional loci of $\phi_{c_5}^{+}$ and
$\phi_{c_5}^{-}$ having dimensions $6$ and $2$,
respectively.
\end{enumerate}
\end{theorem}

\begin{proof}
We first show that $(\bP^1\times\bP^1,cC_{A_7})$
is K-unstable for $c>c_5$. We have $A_{(\bP^1\times\bP^1,cC_{A_7})}(F)= 5-8c$. A direct computation gives
\[
\vol\bigl(2\pi^*H-tF\bigr)
\ =\ 
\begin{cases}
8-\tfrac14t^2,
&0\leq t\leq\tfrac{16}{3},\\[2pt]
2(6-t)^2,
&\tfrac{16}{3}\leq t\leq6,
\end{cases}
\]
and hence
\[
S_{(\bP^1\times\bP^1,cC_{A_7})}(F)
\ =\ \frac{34}{9}-\frac{17}{3}c,\qquad
\beta_{(\bP^1\times\bP^1,cC_{A_7})}(F)
\ =\ \frac{11}{9}-\frac73c\ <\ 0
\] for any $c>c_5$.

\begin{lem}\label{de-non A_7}
The pair $(\bP^1\times\bP^1,cC_{A_7})$ admits a special
degeneration to $(\bP(1,2,9),cC'_{A_7})$.
\end{lem}

\begin{proof}
Let $\pi:\sX\to(\bP^1\times\bP^1)\times\bA^1_t$ be the $(1,1,4)$-weighted
blowup in the local coordinates $(t,z,x')$.
Its central fiber is $S\cup X$, where
$S\simeq\bP(1,1,4)$ is the exceptional divisor and
$X$ is the strict transform of the original central fiber.
Let $\tQ\subset X$ be the strict transform of
$Q=V(q)$. We perform the birational transformations
\[
\begin{tikzcd}
&\sX
  \arrow[ld,"\pi"']
  \arrow[rd,"g"]
  \arrow[rr,dashed,"f"]
&&\sX^+
  \arrow[ld,"h"']
  \arrow[rd,"\psi"]
\\
(\bP^1\times\bP^1)\times\bA^1_t&&\sY&&\sZ,
\end{tikzcd}
\]
where $g$ contracts $\tQ$, $h$ is its flip,
and $\psi$ contracts the strict transform $X'$ of $X$
to a point. On the central fibers, these maps give
\[
\begin{tikzcd}
&S\cup X
  \arrow[ld,"\pi"']
  \arrow[rd,"g"]
  \arrow[rr,dashed,"f"]
&&\hat S\cup X'
  \arrow[ld,"h"']
  \arrow[rd,"\psi"]
\\
\bP^1\times\bP^1&&S\cup X'&&S'.
\end{tikzcd}
\] Write $[u':v':s']$ for the coordinates on
$S\simeq\bP(1,1,4)$. Then
$F=S\cap X=V(u')$, and $\tQ$ meets $F$
at $[0:0:1]$, which has type $\tfrac14(1,1)$ on $S$
and type $\tfrac14(1,3)$ on $X$. On $X$, we have
\[
F^2=-\tfrac14,\qquad
\tQ^2=-\tfrac14,\qquad
F\cdot\tQ=\tfrac34.
\]
After contracting $\tQ$, the image $F'$ of $F$ has self-intersection $2$ on $X'$ and
self-intersection $-2$ on $\hat S$. The morphism $h|_{\hat S}:\hat S\to S$ extracts the
flipped curve $\tQ^+$; a toric model identifies it with the $(9,1)$-weighted blowup at
$[0:0:1]$. Contracting $F'$ gives the toric surface
$S'\simeq\bP(1,2,9)$. Tracking the strict transform of
$C_{A_7}\times\bA^1$ through these transformations
gives the limiting curve $C'_{A_7}=V(s^2+u^2v^8)$.
\end{proof}

By construction, this special test configuration is
induced by the divisorial valuation $\ord_F$.
The pair $(\bP(1,2,9),c_5C'_{A_7})$ is a
$\mathbb T$-pair of complexity one under the $\bG_m$-action
\[\lambda\ \colon \ 
[u:v:s]\ \longmapsto\ [u:\lambda v:\lambda^4s].
\]
There are no horizontal divisors, and the vertical
divisors are $V(u),V(v)$ and $V(s-auv^4)$ for $a\in\bC$. At $c=c_5$, their $\beta$-invariants are positive, and
\[
\Fut_{(\bP(1,2,9),c_5C'_{A_7})}(\lambda)
\ =\ \tfrac{11-21c_5}{9}\ =\ 0.
\]
Thus $(\bP(1,2,9),c_5C'_{A_7})$ is K-polystable by \Cref{thm:complexity-one}. Openness of K-semistability implies that $(\bP^1\times\bP^1,c_5C_{A_7})$ is K-semistable, and the interpolation of K-stability then shows that $(\bP^1\times\bP^1,cC_{A_7})$ is K-stable for $0<c<c_5$.

\smallskip
For any curve
\[
C'=V\bigl(s^2+u^2v^8+u^6g_{12}(u,v)\bigr)
\subseteq\bP(1,2,9),
\]
where $g_{12}$ is weighted homogeneous of degree $12$, there is a $\bG_m$-equivariant degeneration of $C'$ to $C'_{A_7}$. Consequently, $(\bP(1,2,9),c_5C')$ is K-semistable. For general $g_{12}$, the pair $(\bP(1,2,9),\tfrac23C')$ is klt, so interpolation of K-stability implies that $(\bP(1,2,9),cC')$ is K-stable for every $c_5<c<\tfrac23$. On the other hand, the $\beta$-invariant associated with $\lambda$ shows that these pairs are K-unstable for $c<c_5$. Thus the $\phi^+_{c_5}$-exceptional locus $E^+_5$ contains the closure of the locus of these general pairs, which has dimension $6$ by a parameter count. The proper transform of the non-separating $A_7$-locus is contained in $E^-_5$ and has dimension $2$. Applying \cite[Theorem 0.2.5]{DH98} to the local VGIT presentation at $[(\bP(1,2,9),c_5C'_{A_7})]$, we deduce that $E_5$ is zero-dimensional near this point. Hence, locally,
\[
    E_5\ =\ \{[(\bP(1,2,9),c_5C'_{A_7})]\},
\]
$E^-_5$ is the proper transform of the closure of the non-separating $A_7$-locus, and $E^+_5$ is the closure of the locus described above. In particular, $\phi_{c_5}$ is a flip locally over $[(\bP(1,2,9),c_5C'_{A_7})]$.
\end{proof}

\subsection{$A_{6}$-singularity}
After a projective transformation, a $(3,3)$-curve $C_{A_6}$
with an $A_6$-singularity at $p\coloneqq(0,0)$ can be
written in affine coordinates as
\[
r\bigl(z^5+r(z^3+2z^2+xz+4x+4z+8)\bigr)
+a_0q^3+a_1q^2(x+2z)+a_2q(x+2z)^2=0,
\]
where
\[
r\coloneqq x+2z-z^2,\qquad q\coloneqq xz+2x+4z,
\]
and $(a_0,a_1,a_2)\in\bA^2\times(\bA^1\setminus\{0\})$. Choosing the jet coordinate
$x'=x+2z-z^2+\tfrac{z^5}{16}$ at $p$, the equation
becomes
\[
8(x')^2+a_2z^7+(\text{higher order terms})=0,
\]
where $(z,x')$ has weight $(2,7)$.
Let
\[
\pi:X\coloneqq\Bl_{(2,7)}(\bP^1\times\bP^1)
\longrightarrow\bP^1\times\bP^1
\]
be the corresponding weighted blowup at $p$, and
let $F$ be its exceptional divisor. Set $c_6=\tfrac{17}{30}$, and let
$C'_{A_6}\subseteq\bP(1,2,9)$ be the curve defined
by $s^2+u^4v^7=0$.  Let $H \in |\cO_{\bP^1 \times \bP^1}(1,1)|$ denote the class of the hyperplane section.

\begin{theorem}\label{The wall corresponds to A6}
Keep the notation introduced above. Then
\begin{enumerate}
\item[\textup{(1)}]
$(\bP^1\times\bP^1,cC_{A_6})$ is K-semistable
\textup{(}resp.\ K-polystable\textup{)} if and only if
$c\leq c_6$ \textup{(}resp.\ $c<c_6$\textup{)};
\item[\textup{(2)}]
the K-polystable degeneration of
$(\bP^1\times\bP^1,c_6C_{A_6})$ is isomorphic to
$(\bP(1,2,9),c_6C'_{A_6})$;
\item[\textup{(3)}]
$\phi_{c_6}\colon
\ove{\fM}^{\K}_4(c_6-\epsilon)\dashrightarrow
\ove{\fM}^{\K}_4(c_6+\epsilon)$, locally over $[(\bP(1,2,9),c_6C'_{A_6})]$, is a flip, with the
exceptional loci of $\phi_{c_6}^{+}$ and
$\phi_{c_6}^{-}$ having dimensions $5$ and $3$,
respectively.
\end{enumerate}
\end{theorem}

\begin{proof}
We first show that $(\bP^1\times\bP^1,cC_{A_6})$
is K-unstable for $c>c_6$. We have
\[
A_{(\bP^1\times\bP^1,cC_{A_6})}(F)=9-14c.
\]
A direct computation gives
\[
\vol\bigl(2\pi^*H-tF\bigr)
=
\begin{cases}
8-\tfrac1{14}t^2,
&0\leq t\leq\tfrac{28}{3},\\[2pt]
\tfrac14(12-t)^2,
&\tfrac{28}{3}\leq t\leq12,
\end{cases}
\]
and hence for $c>c_6$
\[
S_{(\bP^1\times\bP^1,cC_{A_6})}(F)
 \ =\ \tfrac{64}{9}-\tfrac{32}{3}c,\qquad
\beta_{(\bP^1\times\bP^1,cC_{A_6})}(F)
\ =\ \tfrac{17}{9}-\tfrac{10}{3}c\ <\ 0.
\]

\begin{lem}\label{de-A_6}
The pair $(\bP^1\times\bP^1,cC_{A_6})$ admits a special
degeneration to $(\bP(1,2,9),cC'_{A_6})$.
\end{lem}

\begin{proof}
Let $\sX_0=(\bP^1\times\bP^1)\times\bA^1_t$,
and let $\pi:\sX\to\sX_0$ be the $(1,2,7)$-weighted
blowup in the local coordinates $(t,z,x')$.
Its central fiber is $S\cup X$, where
$S\simeq\bP(1,2,7)$ is the exceptional divisor and
$X$ is the strict transform of the original central fiber.
Let $\tQ\subset X$ be the strict transform of
$Q=V(xz+2(x+2z))$. As in the proof of
\Cref{de-non A_7}, we first contract and flip
$\tQ$, obtaining a central fiber
$\hat S\cup X'$, and then contract $X'$ to a point.

Write $[x_1:x_2:x_3]$ for the coordinates on
$S\simeq\bP(1,2,7)$. Then
$F=S\cap X=V(x_1)$, and $\tQ$ meets $F$
at $[0:0:1]$, which has type $\tfrac17(1,2)$
on $S$ and type $\tfrac17(1,5)$ on $X$.
On $X$, we have
\[
F^2=-\tfrac1{14},\qquad
\tQ^2=-\tfrac47,\qquad
F\cdot\tQ=\tfrac37.
\]
After contracting $\tQ$, the image $F'$
of $F$ has self-intersection $\tfrac14$ on $X'$
and self-intersection $-\tfrac14$ on $\hat S$.
The induced morphism $\hat S\to S$ extracts the
flipped curve $\tQ^+$; a toric model identifies it with the $(9,4)$-weighted blowup
at $[0:0:1]$. Contracting $F'$ gives the toric surface
$S'\simeq\bP(1,2,9)$, and tracking the strict transform of
$C_{A_6}\times\bA^1$ through these transformations
gives the limiting curve $C'_{A_6}$.
\end{proof}

By construction, this special test configuration is
induced by the divisorial valuation $\ord_F$.
The pair $(\bP(1,2,9),c_6C'_{A_6})$ is a
$\mathbb T$-pair of complexity one under the action
\[
\lambda \ \colon \ [u:v:s]\ \longmapsto\ [u:\lambda^2v:\lambda^7s],
\]
with respect to which there are no horizontal divisors, and the vertical
divisors are
\[
V(u),\qquad V(v),\qquad V(s),\qquad
V(s^2-au^4v^7),\qquad a\in\bC^*.
\]
At $c=c_6$, their $\beta$-invariants are positive, and
\[
\Fut_{(\bP(1,2,9),c_6C'_{A_6})}(\lambda)
\ =\ \tfrac{17-30c_6}{9}\ =\ 0.
\]
Thus $(\bP(1,2,9),c_6C'_{A_6})$ is K-polystable by \Cref{thm:complexity-one}. Openness of K-semistability implies that
$(\bP^1\times\bP^1,c_6C_{A_6})$ is K-semistable, and the interpolation of K-stability then shows that
$(\bP^1\times\bP^1,cC_{A_6})$ is K-stable for any $0<c<c_6$.

\smallskip
For any curve
\[
C'\ =\ V\bigl(s^2+u^4v^7+u^8g_{10}(u,v)\bigr)\ 
\subseteq \ \bP(1,2,9),
\]
where $g_{10}$ is weighted homogeneous of degree $10$, there is a $\bG_m$-equivariant degeneration of $C'$ to $C'_{A_6}$. Consequently, $(\bP(1,2,9),c_6C')$ is K-semistable. For general $g_{10}$, the pair $(\bP(1,2,9),\tfrac23C')$ is klt, so interpolation of K-stability implies that $(\bP(1,2,9),cC')$ is K-stable for every $c_6<c<\tfrac23$. On the other hand, the $\beta$-invariant associated with $\lambda$ shows that these pairs are K-unstable for $c<c_6$. Thus the $\phi^+_{c_6}$-exceptional locus $E^+_6$ contains the closure of the locus of these general pairs, which has dimension $5$ by a parameter count. The proper transform of the $A_6$-locus is contained in $E^-_6$ and has dimension $3$. Applying \cite[Theorem 0.2.5]{DH98} to the local VGIT presentation at $[(\bP(1,2,9),c_6C'_{A_6})]$, we deduce that $E_6$ is zero-dimensional near this point. Hence, locally,
\[
    E_6\ =\ \{[(\bP(1,2,9),c_6C'_{A_6})]\},
\]
$E^-_6$ is the proper transform of the closure of the $A_6$-locus, and $E^+_6$ is the closure of the locus described above. In particular, $\phi_{c_6}$ is a flip locally over $[(\bP(1,2,9),c_6C'_{A_6})]$.
\end{proof}

To conclude this subsection, we exclude other K-moduli walls $c_i$ such that a pair $(\bP(1,2,9),(c_i+\epsilon)C)$ is parametrized by the exceptional locus $E_i^{+}$. We prove in \Cref{the wall corresponds to A8,The wall corresponds to non sep A7,The wall corresponds to A6} that $(\bP(1,2,9),cC)$ is K-semistable (resp. K-unstable) if \begin{enumerate}
     \item $C=V\big(s^{2} + v^{9} + u^{2}g_{16}(u,v)\big)$ and $c=c_4$ (resp. $c<c_4$);
     \item $C=V\big(s^{2} + u^{2}v^{8} + u^{4}g_{14}(u,v) \big)$ and $c=c_5$ (resp. $c< c_5$);
     \item $C=V(s^{2} + u^{4}v^{7} + u^{6}g_{12}(u,v) \big)$ and $c=c_6$ (resp. $c<c_6$).
\end{enumerate}
Thus it suffices to prove the following.
\begin{lem}\label{lem:no-other-wall-129}
    If $C=V\big(s^{2}+u^{6}g_{12}(u,v)\big)$ or $C=V\big(usf_{8}(u,v)+g_{18}(u,v)\big)$, then $(\bP(1,2,9),cC)$ is K-unstable for every $c\in(0,\tfrac23)$.
\end{lem}

\begin{proof}
For $C=V\big(s^{2}+u^{6}g_{12}(u,v)\big)$, consider the divisorial valuation $v$ of weights $(1,3)$ in the local cover coordinates $(u,s)$ at $[0:1:0]$. Then
\[
    A_{(\bP(1,2,9),cC)}(v)\ \leq\ 4-6c,
\]
and, denoting by $H_{\bP(1,2,9)}$ a class in the linear system $|\cO_{\bP(1,2,9)}(1)|$, 
\[
\vol\bigl(\pi^{*}H_{\bP(1,2,9)}-tv\bigr)
\ =\
\begin{cases}
    \frac1{18}-\frac16t^2,
    &0\leq t\leq\frac13,\\[2pt]
    \frac1{12}(1-t)^2,
    &\frac13\leq t\leq1.
\end{cases}
\]
It follows that
\[
    S_{(\bP(1,2,9),cC)}(v)
    \ =\ 216\left(1-\tfrac32c\right)
    \int_0^\infty
    \vol\bigl(\pi^*\cO_{\bP(1,2,9)}(1)-tv\bigr)\,dt
    \ =\ \tfrac{16}{3}-8c.
\]
Consequently,
\[
    \beta_{(\bP(1,2,9),cC)}(v)
    \ \leq\ 4-6c-\left(\tfrac{16}{3}-8c\right)
    \ =\ -\tfrac43+2c\ <\ 0
\]
for every $0<c<\tfrac23$.

\smallskip

For $C=V\big(usf_{8}(u,v)+g_{18}(u,v)\big)$, let $\pi:Y\to\bP(1,2,9)$ be the weighted extraction at $p=[0:0:1]$ corresponding to the primitive vector $\frac19(1,2)$, and denote its exceptional divisor by $E$. Then
\[
    A_{\bP(1,2,9)}(E)\ =\ \tfrac13,
    \qquad
    \ord_E(C)\ \geq\ 1.
\]
Indeed, on the local cover $\bA^2_{u,v}$ at $p$, every monomial of $uf_8(u,v)$ has $(1,2)$-weight $9$, while every monomial of $g_{18}(u,v)$ has weight $18$. Moreover,
\[
\vol\bigl(\pi^*H_{\bP(1,2,9)}-tE\bigr)
\ =\
\begin{cases}
    \frac1{18}-\frac92t^2,
    &0\leq t\leq\frac19,\\[2pt]
    0,
    &t\geq\frac19.
\end{cases}
\]
It follows that
\[
    S_{(\bP(1,2,9),cC)}(E)
    \ =\ (12-18c)\cdot18
    \int_0^{1/9}\left(\tfrac1{18}-\tfrac92t^2\right)\,dt
    \ =\ \tfrac89-\tfrac43c.
\]
Consequently,
\[
    \beta_{(\bP(1,2,9),cC)}(E)
    \ \leq\ \tfrac13-c-\left(\tfrac89-\tfrac43c\right)
    \ =\ -\tfrac59+\tfrac13c\ <\ 0
\]
for every $0<c<\frac23$, completing the proof.
\end{proof}

\subsection{Geometry of K-moduli wall-crossing}

In this section, we show that there are no other K-moduli walls and prove \Cref{Main Theorem 1.1}.

\begin{proof}[Proof of \Cref{Main Theorem 1.1}]
    By \Cref{new surfaces in K-moduli spaces} and \cite[Proposition 11.7]{ABB+23}, we know that the K-moduli stack $\MM_4^\K(c)$ is smooth, and the K-moduli space $\ofM_4^{\K}(c)$ is normal for any $c\in (0,\frac{2}{3})$. Statement (1) follows from \Cref{thm:first-wall}, and statement (3) follows from \Cref{new surfaces in K-moduli spaces}. The geometry of the K-moduli wall-crossings at $c_1,\ldots,c_6$, including the dimensions of their exceptional loci and centers exhibited in \Cref{table:K-moduli walls}, is described in Theorems~\ref{thm:first-wall}, \ref{Second wall crossing}, \ref{thm:the-wall corrsponds-to-sep-A7}, \ref{Main Theorem of X10}, \ref{the wall corresponds to A8}, \ref{The wall corresponds to non sep A7}, and~\ref{The wall corresponds to A6}, with two distinct wall-crossings occurring at $c_4$; we distinguish their associated loci by using primed notation for one of them.

    It remains to show that there are no other wall crossings. Suppose that $c_*\in(0,\frac{2}{3})$ is a K-moduli wall distinct from $c_1,\ldots,c_6$. By \cite[Proposition 3.18]{ADL24}, there exists a pair $(X,c_*D)$ that is K-semistable, whereas $(X,cD)$ is K-unstable for any $c<c_*$. By \Cref{new surfaces in K-moduli spaces}, $X$ is isomorphic to $\bP(1,1,2)$, $\bP(1,2,9)$, or $X_{10}$. However, \Cref{Main Theorem of X10}(2)--(3) shows that the only K-moduli wall arising from pairs $(X_{10},cD)$ is $c_3$, while \Cref{lem:no-other-wall112,lem:no-other-wall-129} exclude any additional walls arising from pairs on $\bP(1,1,2)$ or $\bP(1,2,9)$. This yields a contradiction and completes the proof.
\end{proof}

\section{Moduli of log Calabi--Yau pairs, K3 surfaces, and del Pezzo surfaces} \label{Calabi-Yau section}

In this section, we apply the wall-crossing results for K-moduli spaces established in the previous sections to other moduli problems, and draw
some consequences.

\subsection{Moduli of boundary polarized log Calabi--Yau pairs}
Let $\ove{\fM}^{\CY}_4$ be the good moduli space of the moduli stack of
boundary polarized Calabi--Yau pairs admitting a $\bQ$-Gorenstein smoothing
to a pair $(\bP^1\times\bP^1,\tfrac{2}{3}C)$, where $C$ is a smooth curve of
bidegree $(3,3)$; see \Cref{def:CY4stack}. By \Cref{thm:bpCYwallcrossing}, there is a natural birational morphism
\[
    \Phi\ :\ \ove{\fM}^{\K}_4 \ \longrightarrow \ \ove{\fM}^{\CY}_4 ,
\]
which is an isomorphism over the complement of three closed points of
$\ove{\fM}^{\CY}_4$; see \cite[Remark~6.15 and the proof of
Theorem~6.16]{ADLW25}. We now describe $\Phi$ explicitly over a neighborhood of each of these three points.

\smallskip

Let $\Sigma_7$, $\Sigma_8$, and $\Sigma_{10}$ be the three strata of the GIT
quotient $\ove{\fM}^{\GIT}_4$ introduced in \Cref{strata-for-GIT}. Combining \Cref{prop:klt-is-stable,interpolation,strata-for-GIT}, we see that the only loci contracted by $\Phi$ are the strict transforms of $\ove{\Sigma}_7$, $\ove{\Sigma}_8$, and
$\ove{\Sigma}_{10}$. In particular, $\Phi$ is a small contraction.

\begin{prop} \label{prop:CY-is-ample-model}
    The moduli space $\ove{\fM}^{\CY}_4$ is the ample model of the Hodge line
    bundle on $\ove{\fM}^{\K}_4$.
\end{prop}

\begin{proof}
    This follows immediately from \cite[Theorem~1.1(3)]{blum24} and the functoriality of the Hodge line bundle.
\end{proof}

\subsubsection{$(3,3)$-curves with separating $A_5$-singularities}

The stratum $\Sigma_7$ in $\ove{\fM}^{\GIT}$ is a single point, representing the curve with two separating $A_5$-singularities: \[C_{2A_5}\ =\  V(x_0x_3-x_1x_2,\,x_0x_2^2+x_1^2x_3).\] This curve consists of three components: two rulings of $\bP^1\times\bP^1$ and a smooth $(1,3)$-curve meeting each ruling at an $A_5$-singularity.

\smallskip

\begin{figure}[h]
    \centering
    \begin{subfigure}[c]{.30\textwidth}    \centering\resizebox{\textwidth}{!}{\begin{tabular}{c}
\begin{tikzpicture}[gren0/.style = {draw, circle,fill=greener!80,scale=.7},gren/.style ={draw, circle, fill=greener!80,scale=.4},blk/.style ={draw, circle, fill=black!,scale=.08},plc/.style ={draw, circle, color=white!100,fill=white!100,scale=0.02},smt/.style ={draw, circle, color=gray!100,fill=gray!100,scale=0.02},lbl/.style ={scale=.2}] 

\draw [-,color=gray] (-3.4,-.2) to [out=15, in=-90] (-2.2,.8) to [out=90, in=180] (-1.2,2.8);
\draw [-,color=gray] (-.2,-2.6) to (-3.4,-0.2);
\draw [-,color=gray] (0.2,2.6) to (-0.2,-2.6);
\draw [-] (-1.2,2.8) to (0.2,2.6);

\draw [-,color=blue] (0.2,2.6) to (-.2,-2.6);
\draw [-,color=black] (-1.2,2.8) to (0.2,2.6);

\filldraw (-.2,-2.6) circle (1pt);
\filldraw (0.2,2.6) circle (1pt);
\filldraw (-1.2,2.8) circle (0.4pt);

\draw [-,color=gray] (3.4,.2) to [out=195, in=90] (2.2,-.8) to [out=-90, in=0] (1.2,-2.8);
\draw [-,color=gray] (.2,2.6) to (3.4,0.2);
\draw [-,color=gray] (-0.2,-2.6) to (0.2,2.6);
\draw[-] (1.2,-2.8) to (-0.2,-2.6);

\draw [-,color=blue] (-0.2,-2.6) to (.2,2.6);
\draw [-,color=black] (1.2,-2.8) to (-0.2,-2.6);

\filldraw (.2,2.6) circle (1pt);
\filldraw (-0.2,-2.6) circle (1pt);
\filldraw (1.2,-2.8) circle (0.4pt);

\node[above left, node font=\small] at (-.2,-2.6) {$\frac{1}{3}(1,1)$};
\node[above right, node font=\small] at (-.2,-2.6) {$\frac{1}{3}(1,2)$};
\node[below left, node font=\small] at (.2,2.6) {$\frac{1}{3}(1,2)$};
\node[below right, node font=\small] at (.2,2.6) {$\frac{1}{3}(1,1)$};
\node[below left, node font=\small] at (0,0) {$l_1$};
\node[above right, node font=\small] at (0,0) {$l_1'$};

\end{tikzpicture}
\end{tabular}}
        \caption{$S_{2A_5}$}\label{fig:S2A5}
    \end{subfigure}%
    \hspace{1em}
    \begin{subfigure}[c]{.64\textwidth}
        \centering
        \begin{minipage}[c]{.30\textwidth}  \centering\resizebox{\textwidth}{!}{\begin{tabular}{c}
\begin{tikzpicture}[gren0/.style = {draw, circle,fill=greener!80,scale=.7},gren/.style ={draw, circle, fill=greener!80,scale=.4},blk/.style ={draw, circle, fill=black!,scale=.08},plc/.style ={draw, circle, color=white!100,fill=white!100,scale=0.02},smt/.style ={draw, circle, color=gray!100,fill=gray!100,scale=0.02},lbl/.style ={scale=.2}] 

\draw [-,color=gray] (-3.4,-.2) to [out=15, in=-90] (-2.2,.8) to [out=90, in=180] (-1.2,2.8);
\draw [-,color=gray] (-.2,-2.6) to (-3.4,-0.2);
\draw [-,color=gray] (0.2,2.6) to (-0.2,-2.6);
\draw [-] (-1.2,2.8) to (0.2,2.6);

\draw [-,color=blue] (0.2,2.6) to (-.2,-2.6);
\draw [-,color=black] (-1.2,2.8) to (0.2,2.6);

\filldraw (-.2,-2.6) circle (1pt);
\filldraw (0.2,2.6) circle (1pt);
\filldraw[black] (-1.2,2.8) circle (0.4pt);

\node[below left, node font=\normalsize] at (-.2,-2.8) {$S_{A_5}$};
\node[right, node font=\small] at (-.2,-2.6) {$\frac{1}{3}(1,1)$};
\node[below right, node font=\small] at (.2,2.6) {$\frac{1}{3}(1,2)$};
\node[above, node font=\small] at (-.5,2.7) {$\pi^{-1}(p_1)$};
\node[below right, node font=\small] at (0,0) {$l_1$};
\end{tikzpicture}
\end{tabular}}
        \end{minipage}%
        \hspace{.1em}
        $\xrightarrow{\text{\scriptsize wt.\ blow-up of $p_1$}}$%
        \hspace{.1em}
        \begin{minipage}[c]{.38\textwidth}
            \centering\resizebox{\textwidth}{!}{\begin{tabular}{c}
\begin{tikzpicture}[gren0/.style = {draw, circle,fill=greener!80,scale=.7},gren/.style ={draw, circle, fill=greener!80,scale=.4},blk/.style ={draw, circle, fill=black!,scale=.08},plc/.style ={draw, circle, color=white!100,fill=white!100,scale=0.02},smt/.style ={draw, circle, color=gray!100,fill=gray!100,scale=0.02},lbl/.style ={scale=.2}] 

\draw [-,color=gray] (-3.4,.2) to [out=0, in=-90] (-1.6,1.4) to [out=90, in=180] (0.2,2.6);

\draw [-,color=gray] (-.2,-2.6) to (-3.4,.2);
\draw [-,color=gray] (0.2,2.6) to (-.2,-2.6);

\draw [-,color=blue] (0.2,2.6) to (-.2,-2.6);

\filldraw (-.2,-2.6) circle (1pt);
\filldraw (.2,2.6) circle (1pt);

\node[below left, node font=\normalsize] at (-.2,-2.8) {$\bP(1,1,3)$};
\node[right, node font=\small] at (-.2,-2.6) {$\frac{1}{3}(1,1)$};
\node[right, node font=\small] at (0.2,2.6) {$p_1$};
\node[right, node font=\small] at (0,0) {$(x_1 = 0)$};
\end{tikzpicture}
\end{tabular}}
        \end{minipage}
        \caption{The construction of $S_{A_5}$}\label{fig:SA5}
    \end{subfigure}
    \caption{The surfaces $S_{2A_5}$ and $S_{A_5}$.}
    \label{fig:A5-surfaces}
\end{figure}
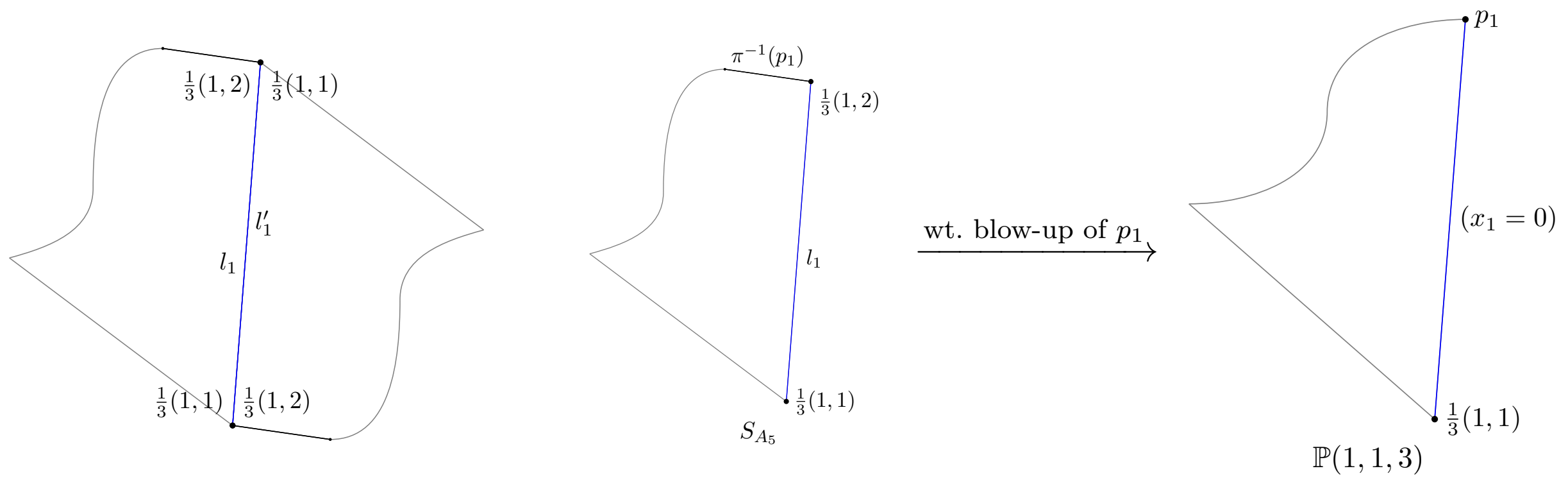

    Let $\pi:S_{A_5}\to \bP(1,1,3)$ be the weighted blowup of $\bP(1,1,3)_{[x_0:x_1:x_2]}$ at $p_1 = [1:0:0]$ with weight $(1,3)$ in $(x_1, x_2)$, depicted in Figure \ref{fig:SA5}. Denote by $\ell_{1}\subset S_{A_5}$ the strict transform of $(x_1=0)\subset \bP(1,1,3)$. Let $[u_0:u_1]$ be the coordinates of $\ell_1$ such that $\pi([u_0:u_1]) = [u_0:0:u_1^3]$. Let $S_{2A_5}$ be the surface obtained by gluing two $S_{A_5}$ surfaces along $\ell_1$ via the automorphism $[u_0:u_1]\mapsto [u_1:u_0]$, depicted in Figure \ref{fig:S2A5}. Let $C_{2A_5}^\circ$ be the curve on $S_{2A_5}$ defined by gluing the two strict transforms of the curves $V(x_2^2 - x_0^3 x_2)$.

\begin{prop}\label{prop:C2A5-replacement}
    We keep the notation introduced above. Then:
    \begin{enumerate}
        \item[\textup{(1)}] the pair $(\bP^1\times \bP^1, (\tfrac{2}{3}-\epsilon)C_{2A_5})$ is K-polystable for any $0<\epsilon\ll1$;
        \item[\textup{(2)}] the polystable degeneration of $(\bP^1\times \bP^1, \tfrac{2}{3}C_{2A_5})$ in $\MM^{\CY}_4$ is $(S_{2A_5}, \tfrac{2}{3}C_{2A_5}^\circ)$; and
        \item[\textup{(3)}] $\Phi$ is an isomorphism over an open neighborhood of $[(S_{2A_5}, \tfrac{2}{3}C^\circ_{2A_5})]$.
    \end{enumerate}
\end{prop}

\begin{proof}
    This is the content of \cite[Lemma~6.1, Propositions~6.3 and~6.7]{ADLW25}.
\end{proof}

\subsubsection{$(3,3)$-curves with $D_4$-singularities}

The stratum $\Sigma_{10}$ in $\ove{\fM}^{\GIT}_4$ parametrizes the unions of
three conics meeting in two $D_4$-singularities; by \Cref{strata-for-GIT}, it
is a two-dimensional family of curves.

\smallskip

    Work in $\bP^3$ with coordinates $[x_0:x_1:x_2:x_3]$. Let
    \[
        S_{2D_4} \ =\ \bP^2_{x_0}\cup\bP^2_{x_3}\ =\ V(x_0x_3)\ \subseteq\ \bP^3
    \]
    be the non-normal quadric surface obtained by gluing two planes along a
    line, and let $C^{\circ}_{2D_4} = V(x_1^3+x_2^3)|_{S_{2D_4}}$ be the curve
    cut out on it by $x_1^3+x_2^3$. For a nonzero quadratic form
    $q = q(x_1,x_2)$, set
    \[
        X_q \ =\ V(x_0x_3+q(x_1,x_2))\ \subseteq\ \bP^3,
        \qquad
        D_q \ =\ V(x_1^3+x_2^3)|_{X_q},
    \]
    so that $X_q$ is a normal quadric surface and $D_q$ has two $D_4$-singularities
    in the smooth locus of $X_q$.

\begin{prop}\label{prop:D4replacement}
    We keep the notation introduced above. Then:
    \begin{enumerate}
        \item[\textup{(1)}] the K-polystable pairs in $\ove{\fM}^{\K}_4$ admitting a special degeneration to $(S_{2D_4},\tfrac{2}{3}C^{\circ}_{2D_4})$ in $\ove{\fM}^{\CY}_4$ are, up to a coordinate change, exactly the pairs $(X_q,\tfrac{2}{3}D_q)$, and they form a two-dimensional locus in $\ove{\fM}^{\K}_4$;
        \item[\textup{(2)}] the pair $(S_{2D_4},\tfrac{2}{3}C^{\circ}_{2D_4})$ is polystable in $\ove{\fM}^{\CY}_4$; and
        \item[\textup{(3)}] $\Phi$ is an isomorphism over a punctured neighborhood of $[(S_{2D_4},\tfrac{2}{3}C^{\circ}_{2D_4})]$, and the fiber of $\Phi$ over $[(S_{2D_4},\tfrac{2}{3}C^{\circ}_{2D_4})]$ is the two-dimensional locus in \textup{(1)}, which is the strict transform of $\ove{\Sigma}_{10}$ under the birational map $\ove{\fM}^{\GIT}_4 \dashrightarrow \ove{\fM}^{\K}_4$.
    \end{enumerate}
\end{prop}

\begin{proof}
    This is the content of \cite[Lemmas~6.9 and~6.10, Proposition~6.12]{ADLW25}.
\end{proof}

\subsubsection{$(3,3)$-curves with a non-separating $A_5$-singularity}

The stratum $\Sigma_8$ in $\ove{\fM}^{\GIT}_4$ parametrizes the $(3,3)$-curves with a non-separating $A_5$-singularity; by \Cref{strata-for-GIT}, it is a
four-dimensional family of curves.

\smallskip

    Let $[a_0:a_1:a_2:a_3]$ be the weighted coordinates on $\bP(1,1,2,3)$, let
    \[
        S_{\mathrm{ns}} \ =\ \bP(1,2,3)\cup\bP(1,1,3)\ =\ V(a_1a_2)
        \ \subseteq\ \bP(1,1,2,3),
    \]
    and let $C^{\circ}_{\mathrm{ns}} = V(a_1a_2,\,a_0^6+a_3^2)$ be the curve on
    $S_{\mathrm{ns}}$.

\begin{prop}\label{prop:A5ns-replacement}
    We keep the notation introduced above. Then:
    \begin{enumerate}
        \item[\textup{(1)}] a point $[(X,\tfrac{2}{3}D)]\in\ove{\fM}^{\K}_4$ admits a weakly special degeneration to $(S_{\mathrm{ns}},\tfrac{2}{3}C^{\circ}_{\mathrm{ns}})$ in $\ove{\fM}^{\CY}_4$ if and only if it lies in the strict transform of $\ove{\Sigma}_8$ under the birational map $\ove{\fM}^{\GIT}_4\dashrightarrow\ove{\fM}^{\K}_4$;
        \item[\textup{(2)}] the pair $(S_{\mathrm{ns}},\tfrac{2}{3}C^{\circ}_{\mathrm{ns}})$ is polystable in $\ove{\fM}^{\CY}_4$; and
        \item[\textup{(3)}] $\Phi$ is an isomorphism over a punctured neighborhood of $[(S_{\mathrm{ns}},\tfrac{2}{3}C^{\circ}_{\mathrm{ns}})]$, and the fiber of $\Phi$ over $[(S_{\mathrm{ns}},\tfrac{2}{3}C^{\circ}_{\mathrm{ns}})]$ is the locus in \textup{(1)}, which has dimension $4$ in $\ove{\fM}^{\K}_4$.
    \end{enumerate}
\end{prop}

\begin{proof}
    Part \textup{(2)} is \cite[Lemma~6.14]{ADLW25}.

    We prove the forward implication in \textup{(1)}. By
    \Cref{new surfaces in K-moduli spaces}, for every
    $[(X,\tfrac{2}{3}D)]\in\ove{\fM}^{\K}_4$ the surface $X$ is isomorphic to
    one of $\bP^1\times\bP^1$, $\bP(1,1,2)$, $X_{10}$, or $\bP(1,2,9)$. If
    $(X,\tfrac{2}{3}D)$ is klt, then it is stable in $\ove{\fM}^{\CY}_4$ by
    \cite[Theorem~6.15]{ABB+23} and admits no nontrivial special degeneration;
    hence $(X,\tfrac{2}{3}D)$ is strictly log canonical. Since $D$ is a
    deformation of $C^{\circ}_{\mathrm{ns}}$, it follows from
    \Cref{prop:klt-is-stable,interpolation,strata-for-GIT,prop:C2A5-replacement,prop:D4replacement}
    that $D$ has a non-separating $A_5$-singularity provided that $X$ is isomorphic to $\bP^1\times\bP^1$. Therefore, it suffices to prove that for the other $X$, the following result holds.
    
\begin{prop}
    The pair $(X,\tfrac{2}{3}D)$ deforms to $(\bP^1\times\bP^1,\tfrac{2}{3}C)$ for some $C$ with a non-separating $A_5$-singularity.
\end{prop}

\begin{proof}We first show that every log canonical center of
$(X,\tfrac23D)$ lies in the smooth locus of $X$.
If $X\simeq\bP(1,1,2)$, then
\Cref{curves meet the vertex point in a cusp},
together with the local equations on the index-one
cover, shows that the pair is klt at the vertex.
If $X\simeq X_{10}$, then $D$ avoids its unique
singular point by \Cref{lem:curves-on-X10}.
If $X\simeq\bP(1,2,9)$, the classification preceding
\Cref{lem:no-other-wall-129} shows that $D$ avoids
the $\frac19(1,2)$-point, while its pullback to the
local cover at the $\frac12(1,1)$-point is either
absent or has at worst an $A_3$-singularity.
Their log canonical thresholds are greater than
$\frac23$, so the pair is klt there as well.

\smallskip

Let $p$ be a log canonical center of
$(X,\tfrac23D)$. Then $p\in X_{\mathrm{sm}}$,
and $D$ is reduced by log canonicity.
The ordinary blowup at $p$ gives
$\mult_pD\leq3$, while $\mult_pD\geq2$ since the pair is not klt at $p$.
If $\mult_pD=2$, the germ is of type $A_n$, and $\lct_p(X;D)=\frac12+\frac1{n+1}=\frac23$ forces $n=5$. If $\mult_pD=3$ and the tangent cone has a repeated factor, choose coordinates in which its cubic part is divisible by $y^2$. The valuation $v$ of weights $(2,3)$ on $(x,y)$ then satisfies $A_{(X,\frac23D)}(v)\leq5-\frac23\cdot8<0$, contradicting log canonicity. Thus $D$ has an ordinary triple point at $p$, which is of type $D_4$. However, since $C^{\circ}_{\mathrm{ns}}$ has at worst double points, the upper semicontinuity of multiplicity excludes the $D_4$-singularity, and hence $D$ has an $A_5$-singularity at $p$.
\begin{lemma}\label{lem:A5-deformation}
Let $L=\cO_X(D)$, and let $Z$ be the Tjurina subscheme of $D$ supported at $p$. Then
\[
    H^0(X,L)\ \longrightarrow\ H^0(Z,L|_Z)
\]
is surjective.
\end{lemma}

\begin{proof}
Since $(D,p)$ is an $A_5$-singularity, its Tjurina algebra is abstractly isomorphic to $k\llbracket t\rrbracket/(t^5)$. If $X=\bP(1,2,9)$, one can write $D=V(s^2+G_{18}(u,v))$, and $p$ lies in the chart $u=1$. The functions $v^i$ for $0\leq i\leq4$ restrict to a basis of the Tjurina algebra and are restrictions of global sections of $L$. If $X=X_{10}$, an automorphism puts $D$ in the form $V_{X_{10}}(w-g_9(x,y))$. Its $A_5$-point lies in the chart $x=1$, where eliminating $w$ gives $z^2=g_9(1,y)-y^{10}$. The functions $y^i$ for $0\leq i\leq4$ likewise restrict to a basis of the Tjurina algebra and come from global sections of $L$.

Finally, let $X=\bP(1,1,2)_{[u,v,s]}$. We may assume that $p=[1:0:0]$. We claim that the ruling $R=V(v)$ cannot be tangent to $D$ at $p$. Indeed, if $R$ is not a component of $D$, then $(D\cdot R)_p\geq4>3=(D\cdot R)$; if it is a component of $D$, then $((D-R)\cdot R)_p=3>\tfrac52=(D-R)\cdot R$. Thus $v$ generates the maximal ideal of the Tjurina algebra, so $v^i$ for $0\leq i\leq4$ restrict to a basis and come from global sections of $L$.
\end{proof}

Choose a $\bQ$-Gorenstein smoothing $\pi\colon\cX\to(0\in B)$ of $X$ to $\bP^1\times\bP^1$ over a pointed smooth curve with $L$ deforming to a line bundle $\cL$. Moreover, the
relative linear system
\[
    q\colon V=\bP_B(\pi_*\cL)\longrightarrow B
\]
is a projective bundle. Let $v_0\in V$ be the point representing $(X,D)$, and
$\kappa\colon(V,v_0)\to \Def(D,p)$ be the forgetful map. Then \Cref{lem:A5-deformation} implies that
\[d\kappa|_{T_{v_0}V_0}\ \colon\  T_{v_0}V_0\ \longrightarrow \ T_0\Def(D,p)\]
is surjective. Since $q$ is smooth and $\Def(D,p)$ is
smooth, the induced map
\[(q,\kappa)\ \colon\ (V,v_0)\ \longrightarrow\ (B,0)\times \Def(D,p)\]
is smooth. Consequently, $\kappa^{-1}(0)\to B$ is smooth near $v_0$ and hence admits a section through $v_0$. This gives a deformation $D_t$ of $D$ along the smoothing of $X$ in which the $A_5$-germ remains analytically locally trivial. By \Cref{prop:C2A5-replacement}, the $A_5$-singularity of $D_t$ is non-separating. 
\end{proof}

Let $X=\bP^1\times\bP^1$ with coordinates $([x:y],[z:w])$, and $D\subseteq X$ be a $(3,3)$-curve with a non-separating $A_5$-singularity. Up to a change of coordinates we may assume that the singularity is the point $p=([0:1],[0:1])$, and that in the affine chart $y=w=1$ the curve $D$ is given by
    \begin{align*}
        &(x+2z-z^{2})\bigl(z^{5}+(x+2z-z^{2})(z^{3}+2z^{2}+xz+4x+4z+8)\bigr)
        + a_{0}\,(xz+2(x+2z))^{3}\\
        &\quad + a_{1}\,(xz+2(x+2z))^{2}(x+2z)
        + a_{2}\,(xz+2(x+2z))(x+2z)^{2}
        + a_{3}\,(x+2z)^{3} \ =\ 0 ,
    \end{align*}
    where $(a_0,a_1,a_2,a_3)\in\bA^3\times(\bA^1\setminus\{0\})$;  $p$ is the singular
    point of $D$ for every such parameter. A direct computation shows this is indeed a $(3,3)$-form. Choosing the jet coordinate
    $x'=x+2z-z^{2}+\tfrac{z^{5}}{16}$ at $p$, the equation of $D$ becomes
    \[
        8(x')^{2} + a_{3}z^{6} + (\text{higher order terms})\ =\ 0,
    \]
    where $(x',z)$ has weight $(3,1)$. Let
    $Q=V(xz+2(xw+2yz))\subseteq\bP^1\times\bP^1$, a curve of bidegree $(1,1)$
    through $p$, which has weight $3$ with respect to this weighting.

    Let $(\cX,\cD)=(X,D)\times\bA^1_t$ be the trivial family of pairs over
    $\bA^1$, and let $\cX'\to\cX$ be the blow-up of $\cX$ along $Q\times\{0\}$.
    Then $\cX'\to\bA^1$ is a fibration whose general fiber is isomorphic to
    $\bP^1\times\bP^1$ and whose central fiber is $\hat{X}_1\cup\hat{X}_2$,
    where $\hat{X}_1\simeq\bP^1\times\bP^1$ is the strict transform of
    $\bP^1\times\bP^1\times\{0\}$ and $\hat{X}_2$ is isomorphic to the
    Hirzebruch surface $\bF_2$. Let $C$ be a smooth curve in
    $|-K_{\bP^1\times\bP^1}|$ and let $\cC$ be the strict transform on $\cX'$ of
    $C\times\bA^1\subseteq\cX$. Running the minimal model program for
    $(\cX',\cC+\epsilon\hat{X}_1)$ over $\bA^1$, with $0<\epsilon\ll1$ rational,
    contracts the component $\hat{X}_1$ entirely; following $\hat{X}_2$ through
    this procedure contracts its unique $(-2)$-curve, so that the central fiber
    becomes $\bP(1,1,2)_{[u:v:s]}$. Tracking the boundary divisor, the limit
    $D'$ of $D$ is a degree-$6$ curve on $\bP(1,1,2)$ with a
    non-separating $A_5$-singularity in the smooth locus; denote the resulting
    family over $\bA^1$ by $(\cX'',\cD')$.

        The central fiber of $(\cX'',\cD')$ is a degree-$6$ curve on $\bP(1,1,2)$
    with a non-separating $A_5$-singularity at a smooth point; the same
    computation as in the proof of \cite[Lemma~6.14]{ADLW25} then produces a
    weakly special degeneration to $(S_{\mathrm{ns}},\tfrac{2}{3}C^{\circ}_{\mathrm{ns}})$.

    Finally, \textup{(3)} follows from \textup{(1)} and \textup{(2)} together
    with the dimension count of \Cref{strata-for-GIT}.
\end{proof}

Combining \Cref{prop:C2A5-replacement,prop:D4replacement,prop:A5ns-replacement},
we obtain the following \Cref{table:K-to-CY-wall-crossing} on the description of $\Phi$.

\begin{table}[h]
    \renewcommand{\arraystretch}{1.4}
    \centering
    \begin{tabular}{|c|c|c|c|}
        \hline
        \textbf{stratum} & \textbf{dimension} & \textbf{polystable pair in $\ove{\fM}^{\CY}_4$} & \textbf{behavior of $\Phi$}\\
        \hline
        $\Sigma_7$ & $0$ & $(S_{2A_5},\tfrac{2}{3}C^{\circ}_{2A_5})$ & isomorphism\\
        \hline
        $\ove{\Sigma}_8$ & $4$ & $(S_{\mathrm{ns}},\tfrac{2}{3}C^{\circ}_{\mathrm{ns}})$ & small contraction\\
        \hline
        $\ove{\Sigma}_{10}$ & $2$ & $(S_{2D_4},\tfrac{2}{3}C^{\circ}_{2D_4})$ & small contraction\\
        \hline
    \end{tabular}
    \caption{The wall crossing $\Phi\colon\ove{\fM}^{\K}_4\to\ove{\fM}^{\CY}_4$.}
    \label{table:K-to-CY-wall-crossing}
\end{table}

\subsection{Moduli of K3 surfaces}

In this section, we study the Hassett--Keel--Looijenga program for the moduli
space of K3 surfaces arising as cyclic triple covers of $\bP^1\times\bP^1$
branched along $(3,3)$-curves, and we prove \Cref{thm:intro-K3s}. Let $\mtf{F}$ be the coarse moduli space of
these K3 surfaces, which is a ball quotient by \cite{Kon02}, and let
$\ove{\mtf{F}}^{\bb}$ be its Baily--Borel compactification. Let $H_n$, $H_h$,
and $H_v$ denote the nodal, hyperelliptic, and Gieseker-Petri divisors in $\mtf{F}$.  Here, the Gieseker-Petri divisor parametrizes the K3 surfaces arising as covers of $\bP(1,1,2)$, and is called the \textit{vanishing theta locus} in \cite{CMJL2012}.  By \cite[Theorem 6.16]{ADLW25}, the moduli space $\ove{\fM}^{\CY}_4$ is normal, and it follows that there is a natural isomorphism
\[
    \ove{\mtf{F}}^{\bb}\ \cong\ \ove{\fM}^{\CY}_4 .
\]
We use it to identify these two moduli spaces, and accordingly regard $\Phi$ as a morphism $\Phi\colon\ove{\fM}^{\K}_4\to\ove{\mtf{F}}^{\bb}$. We denote by
$\lambda_{\Hodge}$ the Hodge line bundle on $\ove{\mtf{F}}^{\bb}$.

\begin{lemma}[Borcherds' relation] \label{Borcherd relation} On $\ove{\mtf{F}}^{\bb}$, one has the relation \begin{equation}\label{eq:Borcherds-relation}
    102\lambda_{\Hodge}\ =\ 2H_n + 9H_v + 28H_h.
\end{equation}
\end{lemma}

\begin{proof}
    By \cite[Theorem 5.11]{CMJL2012}, the divisor associated to an automorphic form $\ove{\phi}$ on $\ove{\mtf{F}}^{\bb}$ is $H_n + \tfrac{9}{2}H_v + 14H_h$. By \cite[Proof of Theorem 1.2]{Borcherds1997FamiliesOK}, $\divisor(\ove{\phi})=a_{\Lambda}\lambda_{\Hodge}$, where $a_{\Lambda}=12+\tfrac{\#R(\Lambda^\perp)}{2}$. By \cite[Lemma 5.4]{CMJL2012}, $\Lambda^\perp$ is isometric to $A_2\oplus E_6$, which has $78$ roots; hence the statement follows. 
\end{proof}

\begin{theorem}[Hassett--Keel--Looijenga program for $\ove{\mtf{F}}^{\bb}$] \label{HKL program}
    For any $c \in (0, \tfrac{2}{3}) \cap \bQ$, set $t(c)\coloneqq\frac{2-3c}{72c}$. Then the CM line bundle on $\ove{\fM}_4^{\K}(c)$ is proportional to $\lambda_c+t(c)(4H_{v,c} + \tfrac{32}{3}H_{h,c})$. In particular, one has
\[ \ove{\fM}_4^{\K}(c) \ \simeq \ \Proj_{\bC}R\big(\mtf{F},\ \lambda_{\Hodge}+t(c)(4H_{v} + \tfrac{32}{3}H_{h})\big).\]
\end{theorem}

To prove the theorem, we first make some preparation by computing certain intersection numbers. Let $C$ be a general hyperelliptic curve of genus four. Then we can realize $C$ as a curve
\[
C=V(z^2-f_{10}(x,y))\ \subseteq \ S \coloneqq \ \bP(1,1,5),
\]
where $f_{10}$ is a binary form with ten distinct roots. Let $B\coloneqq |g^1_2|\simeq\bP^1$ be the hyperelliptic linear series. For each $t\in B$, denote by $F_t$ the corresponding
ruling of $S$, and let $Z_t=C\cap F_t$ be the corresponding divisor. Blowing up $S$ along $Z_t$ and subsequently contracting the proper transform of $F_t$ produces a surface $S_t$ containing a curve $C_t\simeq C$. If $Z_t$ consists of two distinct points, then $S_t\simeq X_{10}$. If $Z_t=2p$ for a Weierstrass point $p\in C$, then $S_t\simeq\mathbb P(1,2,9)$. Thus this construction defines a family of pairs
\[
f:(\cX,\cD)\ \longrightarrow\ B,
\]
and hence a test curve $\Gamma\subset H_h$ parametrizing the different embeddings of the fixed
abstract curve $C$ (cf. \cite[Theorem 2]{Kon02}).

\begin{lemma}\label{lem:intersection}
Let $\Gamma\subseteq H_h$ be the test curve constructed
above. Then
\[
(\Gamma\cdot\lambda_{\Hodge})=\tfrac13,
\qquad
(\Gamma\cdot H_n)=0,
\qquad
(\Gamma\cdot H_v)=10,
\qquad
(\Gamma\cdot H_h)=-2.
\]
\end{lemma}

\begin{proof}
We first compute the Hodge degree. Let $\cR\subseteq S\times B$ be the universal ruling, and let $\cC\coloneqq C\times B$. Write $H$ and $T$ for pullback classes in the linear systems $|\cO_S(1)|$ and $|\cO_B(1)|$ on $S\times B$ respectively. Then
\[
\cR\ \sim\ H+T,
\qquad
\cC\ \sim\ 10H.
\]
Let $\cZ: = {\rm Red}(\cC\cap\cR)$, being naturally isomorphic to $C$.
Consider the blowup
\[
\pi:\cY=\Bl_{\cZ}
(S\times B)\ \longrightarrow\  S\times B
\]
with exceptional divisor $E$, and denote by
$\widetilde{\cC}$ and $\widetilde{\cR}$ the strict transforms. Since $K_S\sim-7H$, we have
\[
K_{\cY/B}\ \sim \ -7 \pi^{*}H+E,
\qquad
\widetilde{\cC}\ \sim \ 10 \pi^{*}H-E.
\]
Consequently,
\[
K_{\cY/B}
 +\tfrac23\widetilde{\cC}
\ \sim_{\bQ}\ -\tfrac13 \pi^{*}H+\tfrac13E \ \sim_{\bQ} \ -\tfrac13\widetilde{\cR}
  +\tfrac13 \pi^{*}T,
\]
where the second equality follows from
$\widetilde{\cR}\sim \pi^{*}(H+T)-E$. The divisor
$\widetilde{\cC}$ is relatively base-point free, and induces a contraction of $\widetilde{\cR}$. Let $f:\cX\to B$ be the induced family, which satisfies
\[
\cO(3K_{\cX/B}+2\cD)
\ \sim_{\bQ} f^*\cO_B(1),
\] by the previous computation and hence $(\Gamma\cdot\lambda_{\Hodge})=\tfrac13$.

Next, we compute the intersection with $H_v$.
The surface $X_{10}$ has a unique cyclic quotient
singularity of type $\frac19(1,2)$, whereas
$\mathbb P(1,2,9)$ has an additional $A_1$-singularity. Thus by Riemann--Hurwitz there are ten (distinct) points on $\Gamma$ at which $Z_t$ becomes nonreduced, corresponding to the intersection of $\Gamma$ with $H_v$. We need to show that these intersections are transverse. Since $\widetilde{\cR}$ is disjoint with the ten $A_1$-singularities on these fibers, $\cX$ is smooth at these $A_1$-singularities. Under the local identification of $H_v$ with
the discriminant divisor parametrizing singular surfaces with an $A_1$-singularity, each branch
point contributes intersection multiplicity one.
Since $C\to\mathbb P^1$ has ten distinct branch
points, it follows that $(\Gamma\cdot H_v)=10$. 

Since the abstract curve $C$ remains smooth throughout the family, $\Gamma$ is disjoint from the nodal divisor $H_n$, and hence $(\Gamma\cdot H_n)=0$. Finally, intersecting the Borcherds relation \Cref{eq:Borcherds-relation} with $\Gamma$ gives $(\Gamma\cdot H_h)=-2$.
\end{proof}

For every $c\in(0,\tfrac{2}{3})\cap\bQ$, let
$\phi_c\colon \ove{\fM}^{\CY}_4 \dashrightarrow \ove{\fM}^{\K}_4(c)$ be the
natural birational contraction, and set
\[\lambda_c\coloneqq (\phi_c)_*\lambda_{\Hodge}, \qquad H_{n,c}\coloneqq (\phi_c)_*H_n, \qquad H_{v,c}\coloneqq (\phi_c)_*H_v, \qquad H_{h,c}\coloneqq (\phi_c)_*H_h.\] By construction,
$\phi_c$ contracts $H_h$ if and only if
$c\in(0,\tfrac{16}{33}]\cap\bQ$, and it contracts $H_v$ if and only if
$c\in(0,\tfrac{1}{6}]\cap\bQ$.

\begin{proof}[Proof of \Cref{HKL program}]
From the construction, for all $0<c<c'<\tfrac{2}{3}$ the birational map $\ove{\fM}_4^{\K}(c')\dashrightarrow \ove{\fM}_4^{\K}(c)$ is a birational contraction, i.e. it does not extract any divisors. Thus by the same proof as in \cite[Theorem 9.4]{ADL24}, it suffices to show that the pullback $(\Phi^{-1})^{*}\Lambda_{\frac{2}{3}-\epsilon,c}$ on $\ove{\mtf{F}}^{\bb}$ is proportional to
\begin{equation*}
    (1-\tfrac{3}{2}c)(4H_{v} + \tfrac{32}{3}H_{h}) +36c\lambda_{\Hodge}.
\end{equation*}
By \cite[Proposition 3.35]{ADL24}, \begin{equation}\label{eq:CM-line-bundle}
        \textstyle (1-\tfrac{3}{2}c)^{-2}\Lambda_{\frac{2}{3}-\epsilon,c} \ = \ (1-\tfrac{3}{2}c)\Lambda_{\frac{2}{3}-\epsilon,0}+36c\Lambda_{c,\Hodge}.
    \end{equation} 
By the cyclic cover construction, one has $(\Phi^{-1})^*\Lambda_{c,\Hodge}=\lambda_{\Hodge}$. As $\Phi^{-1}\big(\mtf{F}\setminus(H_v \cup H_h)\big)$ parametrizes pairs $(\bP^1\times\bP^1, C)$, the pullback of $\Lambda_{\frac{2}{3}-\epsilon,0}$ to $\mtf{F}\setminus(H_v \cup H_h)$ is trivial, since the underlying family of del Pezzo surfaces is an isotrivial $\bP^1\times\bP^1$-fibration. Thus the support of $(\Phi^{-1})^{*}\Lambda_{\frac{2}{3}-\epsilon,0}$ is contained in $H_v\cup H_h$, and we may write $(\Phi^{-1})^{*}\Lambda_{\frac{2}{3}-\epsilon,0}= a_vH_v+a_hH_h$ for some $a_v,a_h\in \bQ$.

Let $\ove{\fM}_4^{\Chow}$ be the GIT quotient of the Chow variety of the canonically embedded genus four curves (cf. \cite{CMJL2012}), and $q: \ove{\fM}_4^{\Chow} \dashrightarrow \ove{\mtf{F}}^{\bb}$ be the natural birational period map. Then by \cite[Theorem 6.2]{CMJL2012}, one can resolve $q$ by blowing up the point $O$ representing the orbit of the ribbon, and contracting the strict transform of a rational curve containing $O$. It follows from \eqref{eq:Borcherds-relation} and \cite[Proposition 1.1, Theorem 2.7 and Theorem 6.2]{CMJL2012}, the natural polarization on $\ove{\fM}_4^{\Chow}$ is proportional to $(q^{-1})_*\lambda_{\Hodge}$.  By \cite[Theorem 6.2]{CMJL2012} again, the exceptional divisor of $q^{-1}$ is $H_h$. Hence one can find $b_h \in \bQ$ such that
\begin{equation} \label{b_h}
    (\lambda_{\Hodge} + b_h H_h)|_{H_h} \ \sim_{\bQ}\  0.
\end{equation} 
Intersecting the curve $\Gamma\subseteq H_h$ constructed in \Cref{lem:intersection} with both sides of \eqref{b_h}, we obtain $b_h = \tfrac{1}{6}$. By the description of the K-moduli wall crossing, the birational map $H_h\dashrightarrow H_{h,c}$ is a birational contraction, and hence
\begin{equation}\label{eq:linear-relation}
    (\lambda_c + \tfrac{1}{6}H_{h,c})|_{H_{h,c}} \ \sim_{\bQ} \ 0.
\end{equation}
By the computation in \Cref{Main Theorem of X10}, the K-moduli wall at which the proper transform of $H_{h}$ in $\ove{\fM}_4^{\K}(c)$ appears first is $c_3=\tfrac{16}{33}$. Moreover, by the geometry of K-moduli wall crossing (cf. \Cref{table:K-moduli walls,strata-for-GIT}), $H_{h,c_3+\epsilon}$ and $H_{v,c_3+\epsilon}$ are disjoint in $\ove{\fM}_4^{\K}(c_3+\epsilon)$. Since the CM line bundle on $H_{h,c_3}$ is $\bQ$-linearly trivial, \Cref{eq:CM-line-bundle,eq:linear-relation} together imply that 
\begin{equation*}
    \frac{1-\tfrac{3}{2}c_3}{36c_3}\,a_h \ =\  \tfrac{1}{6},
\end{equation*}
and hence $a_h = \tfrac{32}{3}$. Similarly, consider the family of curves obtained by intersecting a general pencil of quadrics in $\mathbb P^3$ with a fixed smooth cubic surface. Computing the degree of the CM line bundle of the corresponding family of quadric surfaces, and using the fact that the discriminant divisor of singular quadrics has degree $4$, we obtain $a_v=4$. This finishes the proof.
\end{proof}

\subsection{K-moduli space of del Pezzo surfaces of degree 1} \label{subsection: K-moduli del pezzo surface of degree 1}

In this subsection, we use wall-crossing for the K-moduli of pairs to give an explicit description of the K-moduli space $\fM^{\K}_{\textup{dP}^1}$ of degree one del Pezzo surfaces, which was constructed and studied analytically in \cite{OSS16}.  We first define a GIT moduli space of curves on $\bP(1,1,2)$.

\smallskip
Let $\cU\subseteq\left|\cO_{\mathbb P(1,1,2)}(6)\right|$
be the open locus parametrizing curves that do not pass through the cone point
and are not a triple conic. Every curve in $\cU$ is defined by an equation
\[
as^3+s^2f_2(u,v)+sf_4(u,v)+f_6(u,v)\ =\ 0,
\]
where $a\neq 0$ and $f_i$ is a binary form of degree $i$. After rescaling the
equation, we may assume that $a=1$. The automorphism group of
$\mathbb P(1,1,2)$ is not reductive, and its unipotent radical $N\coloneqq H^0\left(\mathbb P^1,\cO_{\mathbb P^1}(2)\right)$ acts by $s\mapsto s+q_2(u,v)$. The action of $N$ on $\cU$ is free, and every $N$-orbit meets the locus $f_2=0$ in exactly one point: the change of coordinates
\[
s\ \longmapsto \ s-\tfrac{1}{3}f_2(u,v)
\]
puts the equation uniquely into the Weierstrass form
\[
s^3+sg_4(u,v)+g_6(u,v)\ =\ 0.
\]
Consequently, taking the quotient by the unipotent radical identifies the
resulting parameter space with $W\setminus\{0\}$, where
\[
W\coloneqq
\Sym^4 V
\oplus
\Sym^6 V,
\qquad
V\coloneqq H^0\bigl(\bP^1,\cO_{\bP^1}(1)\bigr).
\]
The residual reductive group is $\GL(V)/\boldsymbol{\mu}_2$. Thus we consider the action of $\SL(V)\times\mathbb G_m$ on $W$ by
\[
(g,t)\cdot(f_4,f_6)
\ \coloneqq\ 
\left(
t^4 f_4\circ g^{-1},
t^6 f_6\circ g^{-1}
\right).
\] Let $
\cP\coloneqq\left[
\bigl(W\setminus\{0\}\bigr)/\bG_m
\right]$ be the associated weighted projective stack. As the $\SL(V)$-action on $W$ commutes with the $\bG_m$-action, it descends to an action on $\bP(2^5,3^7)$. We define
\[
\fN^{\GIT}
\ \coloneqq\ 
\bP(2^5,3^7)\sslash_{\cO_{\bP}(1)}\SL(V),
\] which is a normal projective variety.

\begin{defn}
For any $0<c<\frac{2}{3}$, let $\fN^{\K}(c)$ denote the irreducible closed subscheme of $\ove{\fM}_4^{\K}(c)$, endowed with its reduced scheme structure,
whose general point parametrizes a pair with underlying surface $\bP(1,1,2)$.
\end{defn}
\begin{prop}\label{first wall in degree 1 del pezzo}
For any $c\in(\tfrac{1}{6},\tfrac{1}{3})$, there is a natural isomorphism $
\nu_c\colon \fN^{\GIT}\longrightarrow\fN^{\K}(c)$
\end{prop}

\begin{proof}
This follows from the same proof as \cite[Section 5.2]{ADL24}.
\end{proof}

\begin{prop} \label{iso for double cover}
The K-moduli space $\fM^{\K}_{\mathrm{dP}^1}$ is naturally isomorphic
to $\fN^{\K}\bigl(\frac{1}{2}\bigr)$.
\end{prop}

\begin{proof}
The same argument as in \cite[Proposition 6.12]{ADL22} yields a natural
bijective morphism $
\iota\colon
\fN^{\K}\left(\frac{1}{2}\right)
\rightarrow
\fM^{\K}_{\mathrm{dP}^1}$. Since klt surface singularities are quotient singularities, $\bQ$-Gorenstein deformations of klt Fano surfaces are unobstructed by \cite[Proposition 3.1]{HP10}. Hence the K-moduli stack
$\MM^{\K}_{\mathrm{dP}^1}$ is smooth, and its good moduli space
$\fM^{\K}_{\mathrm{dP}^1}$ is normal. Zariski's Main Theorem then
implies that $\iota$ is an isomorphism.
\end{proof}

\begin{proof}[Proof of \Cref{thm:intro-dp1}]
    The diagram follows from Theorem \ref{Main Theorem 1.1}, Proposition \ref{first wall in degree 1 del pezzo} and Proposition \ref{iso for double cover}. In particular, by Table \ref{table:K-moduli walls}, the divisorial contraction $\phi$ is the restriction of the K-moduli wall crossing map decreasing from $c = \frac{1}{3}+\epsilon$ to $c = \frac{1}{3}$ to ${\fN^{\K}(c)}$, and the small contraction $\psi$ is the restriction of the K-moduli wall crossing map increasing from $c = \frac{1}{2}-\epsilon$ to $c = \frac{1}{2}$ to ${\fN^{\K}(c)}$.
\end{proof}

\bibliographystyle{amsalpha}
\bibliography{citation}

\end{document}